\documentclass[a4paper,12pt,frenchb]{article}

\usepackage{amsmath,amsbsy,amsfonts,amssymb,amsthm}
\usepackage[french]{babel}
\newtheorem*{thm}{Th\'eor\`eme}
\newtheorem*{lem}{Lemme}
\newtheorem*{prop}{Proposition}

\begin{document}

 \title{ Sur les donn\'ees endoscopiques dans le cas de l'endoscopie tordue}
\author{J.-L. Waldspurger}
\date{24 f\'evrier 2022}
 \maketitle
 
 {\bf Abstract.} We give a simple combinatorial description of the elliptic endoscopic data of a twisted space under a group $G$, assuming that $G$ is semi-simple and simply connected.  Assuming the same hypothesis  and that the base field  is a number field, we prove that, if two elliptic endoscopic data are equivalent almost everywhere, then they are equivalent. 
\bigskip

{\bf Introduction}

Consid\'erons un corps de nombres $k$ et un groupe r\'eductif connexe $G$ d\'efini sur $k$. Notons ${\mathbb A}$ l'anneau des ad\`eles de $k$ et consid\'erons de plus un caract\`ere continu $\chi$ de $G({\mathbb  A})$ trivial sur $G(k)$. On sait d\'efinir la notion de donn\'ee endoscopique pour $(G,\chi)$. Pour toute place $v$ de $k$, des objets pr\'ec\'edents se d\'eduisent des objets d\'efinis sur le compl\'et\'e $k_{v}$.  Consid\'erons deux donn\'ees endoscopiques sur $k$ et supposons que, pour toute place $v$, les donn\'ees sur $v$ qui s'en d\'eduisent sont \'equivalentes.  B. Lemaire et l'auteur ont prouv\'e en  \cite{LW}, proposition 1, que les deux donn\'ees  \'etaient  alors globalement \'equivalentes. 

On consid\`ere ici le cas de l'endoscopie tordue, c'est-\`a-dire que l'on se donne de plus un espace tordu $\tilde{G}$ sous $G$, d\'efini sur $k$ (cf. \cite{MW} I.1.1 pour une d\'efinition pr\'ecise). Il y a encore une notion de donn\'ee endoscopique. Le r\'esultat de \cite{LW} ne s'\'etend pas \`a ce cas. Nous donnerons un contre-exemple dans le paragraphe \ref{exemples}. Nous prouverons toutefois qu'il reste valide si l'on impose les hypoth\`eses suivantes:

- les donn\'ees endoscopiques sont "$s$-unitaires", cf. paragraphe \ref{ordrefini};

- le groupe $G$ est semi-simple et simplement connexe  (cf. paragraphe \ref{equivalencepresquepartout})ou semi-simple et absolument presque simple (cf. paragraphe \ref{presquesimple}). 
 
 Les donn\'ees elliptiques  sont $s$-unitaires. En fait, le r\'esultat est certainement vrai sans cette hypoth\`ese d'unitarit\'e mais il  nous a  paru  inutile de la lever.  
 
L'int\'er\^et de ce r\'esultat est m\'ediocre. Sa d\'emonstration nous semble plus int\'eressante (acc\`es d'optimisme?). Consid\'erons  le cas o\`u $G$ est  semi-simple et simplement connexe. Dans  la proposition 2 de \cite{LW}, o\`u nous consid\'erions le cas non tordu et o\`u l'on supposait de plus $G$ absolument presque simple, nous avons donn\'e une description facile des  classes d'\'equivalence  des donn\'ees endoscopiques elliptiques. Elle reposait sur une construction assez magique de Langlands, cf. \cite{L} pages 708-709, et celui-ci ne consid\`ere que le cas non tordu. Nous \'etendons ici cette description au cas tordu. Expliquons-la. En fait, tout notre article se passe dans le groupe dual $\hat{G}$. Le corps $k$ n'intervient que par une action continue du groupe de Galois absolu $\Gamma_{k}$ sur $\hat{G}$ par automorphismes alg\'ebriques. On a suppos\'e ci-dessus que $k$ \'etait un corps de nombres mais  la description qui suit vaut aussi si $k$ est  un corps local de caract\'eristique nulle. Puisque $G$ est simplement connexe, $\hat{G}$ est adjoint.  De l'espace tordu $\tilde{G}$ se d\'eduit un automorphisme $\theta$ de $\hat{G}$, qui commute \`a l'action galoisienne. Notons $e$ l'ordre de $\theta$ et $\Theta$ le groupe d'automorphismes engendr\'e par $\theta$. 
 Notons $I(\hat{G})$ l'ensemble des composantes simples de $\hat{G}$. Le groupe $\Theta\times \Gamma_{k}$ agit sur $I(\hat{G})$, notons $i(\hat{G},\theta\times \Gamma_{k})$ le nombre  d'orbites.
 
  Introduisons deux corps locaux $F$  et $E$ de caract\'eristique nulle, de caract\'eristique r\'esiduelle $p>e$, tels que $E$ soit une extension galoisienne cyclique totalement ramifi\'ee de $F$ de degr\'e $e$. Introduisons un groupe r\'eductif connexe ${\bf G}$ d\'efini et quasi-d\'eploy\'e sur $F$  v\'erifiant les hypoth\`eses suivantes: 

${\bf G}$ a m\^emes donn\'ees de racines que $\hat{G}$ (en  oubliant les actions galoisiennes); cela entra\^{\i}ne que $\theta$ se transf\`ere en un automorphisme alg\'ebrique encore not\'e $\theta$ de ${\bf G}$;

  ${\bf G}$ est d\'eploy\'e sur $E$ et il existe un g\'en\'erateur du groupe de Galois $\Gamma_{E/F}$ dont l'action ${\bf G}$ induise une action sur les donn\'ees de racines co\"{\i}ncidant avec $\theta$. 

L'action de $\Gamma_{k}$ sur $\hat{G}$ se transf\`ere en une action de $\Gamma_{k}$ sur ${\bf G}$ par automorphismes alg\'ebriques. On note $\sigma\mapsto \sigma_{G}$ cette action. 
Introduisons les immeubles de Bruhat-Tits $Imm_{E}({\bf G})$ et $Imm_{F}({\bf G})$ de ${\bf G}$ sur $E$ et $F$.  Des actions de $\theta$ et $\Gamma_{k}$ sur ${\bf G}$ se d\'eduisent des actions sur l'immeuble $Imm_{E}({\bf G})$. L'immeuble $Imm_{F}({\bf G})$ s'identifie au sous-ensemble de points fixes $Imm_{E}({\bf G})^{\theta}$. On peut fixer un appartement $App_{E}\subset Imm_{E}({\bf G})$ et une alc\^ove $C\subset App_{E}$ qui soient conserv\'es par $\theta$ et  par l'action de $\Gamma_{k}$. L'ensemble ${\bf A}^{nr}=App_{E}^{\theta}$ est un appartement de $Imm_{F}({\bf G})$ et l'ensemble $C^{nr}=C^{\theta}$ est une alc\^ove de cet appartement.  
On sait bien qu'il y a un ensemble fini de "racines"
$\Delta_{a}^{nr}$ et, pour toute $\beta\in \Delta_{a}^{nr}$, une racine affine $\beta^{aff}$ sur ${\bf A}^{nr}$  de sorte que  

$C^{nr}$ soit l'ensemble des $x\in {\bf A}^{nr}$ tels que $\beta^{aff}(x)>0$ pour toute $\beta\in \Delta_{a}^{nr}$;

pour tout $\beta\in \Delta_{a}^{nr}$, l'hyperplan annulateur de $\beta^{aff}$ soit un mur de l'adh\'erence $\bar{C}^{nr}$ de $C^{nr}$.

Bruhat et Tits d\'efinissent le groupe de Weyl affine \'etendu $W^{aff}$. Il agit sur ${\bf A}^{nr}$. On note $\Omega^{nr}$ le sous-groupe des $v\in W^{aff}$ dont l'action conserve $C^{nr}$. C'est un groupe ab\'elien fini. L'action $\sigma\mapsto \sigma_{G}$ de $\Gamma_{k}$ conserve $C^{nr}$. Notons $\underline{Endo}$ l'ensemble des couples $(\omega_{\star},x)$ tels que

$\omega_{\star}:\Gamma_{k}\to \Omega^{nr}$ est une application continue et l'application $\sigma\mapsto \sigma_{\star}:=\omega_{\star}(\sigma)\circ \sigma_{G}$ est un homomorphisme (\`a valeurs dans le groupe d'automorphismes de ${\bf A}^{nr}$);

$x$ est un \'el\'ement de  $\bar{C}^{nr}$ fix\'e par  $\sigma_{\star}$ pour tout $\sigma\in \Gamma_{k}$.

Pour $(\omega_{\star},x)\in \underline{Endo}$, notons $S(x)$ l'ensemble des $\beta\in \Delta_{a}^{nr}$ tels que $\beta^{aff}(x)=0$. C'est un sous-ensemble propre de $\Delta_{a}^{nr}$. De l'action $\sigma\mapsto \sigma_{\star}$, qui conserve $C^{nr}$, se d\'eduit une action sur $\Delta_{a}^{nr}$, not\'ee de la m\^eme fa\c{c}on. Cette action conserve $S(x)$, donc aussi $\Delta_{a}^{nr}-S(x)$. On v\'erifie facilement que le nombre d'orbites de l'action de $\Gamma_{k}$ sur $\Delta^{a}_{nr}-S(x)$ est minor\'e par $i(\hat{G},\theta\times\Gamma_{k})$. 
Disons que $(\omega_{\star},x)$ est elliptique si ce nombre d'orbites est \'egal \`a $i(\hat{G},\theta\times\Gamma_{k})$.

Deux \'el\'ements $(\omega_{\star,1},x_{1})$ et $(\omega_{\star,2},x_{2})$ de $\underline{Endo}$ sont dits \'equivalents si et seulement s'il existe $\omega\in \Omega^{nr}$ tel que $\omega(x_{1})=x_{2}$ et $\omega\sigma_{\star,1}\omega^{-1}=\sigma_{\star,2}$ pour tout $\sigma\in \Gamma_{k}$. On note $Endo$ l'ensemble des classes d'\'equivalences et $Endo_{ell}$ le sous-ensemble des classes elliptiques. 

Notre r\'esultat est qu'il y a une bijection "naturelle" entre l'ensemble des classes d'\'equivalence de donn\'ees endoscopiques  $s$-unitaires et l'ensemble $Endo$ et que cette bijection se restreint en une bijection entre l'ensemble des classes d'\'equivalence de donn\'ees endoscopiques elliptiques et l'ensemble $Endo_{ell}$.  C'est le th\'eor\`eme du paragraphe \ref{laclassification}. Cela fournit une description combinatoire simple des donn\'ees endoscopiques $s$-unitaires ou elliptiques.

\section{Notations}\label{notations}

On note ${\mathbb N}_{>0}$ l'ensemble des entiers strictement positifs. Pour $n\in {\mathbb N}_{>0}$, on note $\mu_{n}({\mathbb C})$ l'ensemble des  nombres complexes $z$ tels que $z^n=1$. On note $\mu({\mathbb C})$ la r\'eunion des $\mu_{n}({\mathbb C})$ quand $n$ d\'ecrit ${\mathbb N}_{>0}$.  On note $U^1$ le groupe des nombres complexes de module $1$.

Soit $X$ un ensemble sur lequel agit un groupe $G$. On note $X^G$ l'ensemble des points fixes pour cette action. Si $G$ est engendr\'e par un unique \'el\'ement $\theta$, on note aussi $X^{\theta}=X^G$. Pour $x\in X$, on note $Z_{G}(x)$ le fixateur de $x$ dans $G$. Pour tout sous-ensemble $Y\subset X$, on note $Norm_{G}(Y)$ le sous-groupe des \'el\'ements de $G$ qui conservent $Y$.
    
Pour tout groupe alg\'ebrique $G$ d\'efini sur un corps $k$ de caract\'eristique nulle, on note $G^0$ sa composante neutre et $Z(G)$ son centre. Pour un tore $T$ d\'efini sur $k$, on note $X^{*}(T)$, resp. $X_{*}(T)$ le groupe des caract\`eres, resp. cocaract\`eres, de $T$ d\'efinis sur une cl\^oture alg\'ebrique de $k$. 

Soit  $k$ un corps  qui est soit un corps local de caract\'eristique nulle, soit un corps de nombres. On en fixe une cl\^oture alg\'ebrique $\bar{k}$, on note $\Gamma_{k}$ le groupe de Galois de $\bar{k}/k$ et $W_{k}$ le groupe de Weil de $k$. Si $k'$ est une extension galoisienne de $k$, que l'on suppose toujours contenue dans $\bar{k}$, on note $\Gamma_{k'/k}$ le groupe de Galois de $k'/k$.  Si $k$ est un  corps de nombres, on note $V_{k}$ son ensemble de places. Pour toute place  $v\in V_{k}$, on note $k_{v}$ le compl\'et\'e de $k$ en $v$ et on identifie $\Gamma_{k_{v}}$ \`a  un sous-groupe de d\'ecomposition de $\Gamma_{k}$. On a un homomorphisme naturel $W_{k_{v}}\to W_{k}$. 

On fixe pour tout l'article un corps $k$ qui est soit un corps local de caract\'eristique nulle, soit un corps de nombres. 

\section{Donn\'ees endoscopiques}\label{donneesendoscopiques}

  On rappelle ici la notion de donn\'ee endoscopique dans le cas tordu, cf. \cite{MW} I.1.5 pour plus de d\'etails. Soit $G$ un groupe r\'eductif connexe d\'efini sur $k$ et $\tilde{G}$ un espace tordu sous $G$, lui aussi d\'efini sur $k$. On note $\hat{G}$ le dual de Langlands de $G$. Soit ${\bf a}$ un  \'el\'ement de $H^1(W_{k},Z(\hat{G}))$ si $k$ est local, de $H^1(W_{k},Z(\hat{G}))/ker^1(W_{k},Z(\hat{G}))$ si $k$ est un corps de nombres, o\`u $ker^1(W_{k},Z(\hat{G}))$ est le noyau de l'homomorphisme
  $$H^1(W_{k},Z(\hat{G}))\to \prod_{v\in V_{k}}H^1(W_{k_{v}},Z(\hat{G})).$$
   L'espace tordu $\tilde{G}$ d\'etermine un automorphisme $\theta_{Z}$ de $Z(G)$ et on suppose que celui-ci est d'ordre fini. Ces donn\'ees \'etant fix\'ees, on sait d\'efinir la notion de donn\'ee endoscopique du triplet $(G,\tilde{G},{\bf a})$.   Le groupe $\hat{G}$ est muni d'une action de $\Gamma_{k}$, que l'on note $\sigma\mapsto \sigma_{G}$, et il existe par d\'efinition une paire de Borel \'epingl\'ee $\hat{{\cal E}}=(\hat{B},\hat{T},(E_{\alpha})_{\alpha\in \Delta})$  de $\hat{G}$ qui est conserv\'ee par cette action. On fixe une telle paire de Borel \'epingl\'ee $\hat{{\cal E}}$. La donn\'ee de $\tilde{G}$ d\'etermine un automorphisme $\theta$ de $\hat{G}$ qui conserve $\hat{{\cal E}}$ et commute \`a l'action galoisienne. L'ensemble $\hat{G}\theta$ est un espace tordu sous $\hat{G}$. Introduisons le $L$-groupe $^LG=\hat{G}\rtimes W_{k}$. Alors l'ensemble $^LG\theta$ est aussi un espace tordu sous $^LG$. 
   
   Consid\'erons un \'el\'ement semi-simple $\tilde{s}\in \hat{G}\theta$. Posons $\hat{G}'=Z_{\hat{G}}(\tilde{s})^0$. Consid\'erons un sous-groupe ${\cal G}'\subset {^LG} $. On suppose ${\cal G}'\cap \hat{G}=\hat{G}'$. On a donc une suite 
$$(1) \qquad 1\to \hat{G}'\to {\cal G}'\to W_{k}\to 1.$$
On suppose que cette suite est exacte et scind\'ee, c'est-\`a-dire que la projection ${\cal G}'\to W_{k}$ admet  une section  qui soit un homomorphisme continu. On suppose qu'il existe un cocycle $a:W_{k}\to Z(\hat{G})$ dont la classe soit ${\bf a}$ dans $H^1(W_{k},Z(\hat{G}))$ si $k$ est local ou dans  $H^1(W_{k},Z(\hat{G}))/ker^1(W_{k},Z(\hat{G}))$ si $k$ est un corps de nombres, de sorte que, pour tout $(g,w)\in {\cal G}'$, on ait l'\'egalit\'e
$$\tilde{s}(g,w)=(a(w)g,w)\tilde{s}$$
dans $^LG\theta$. 
De la suite (1)  se d\'eduit une $L$-action de $\Gamma_{k}$ dans $\hat{G}'$ que l'on note $\sigma\mapsto \sigma_{G'}$. Soit $G'$ un groupe r\'eductif connexe d\'efini et quasi-d\'eploy\'e sur $k$ de sorte que $\hat{G'}$ muni de cette action galoisienne soit un groupe dual de $G'$. A ces conditions, le triplet $(G',{\cal G}',\tilde{s})$  est appel\'e donn\'ee endoscopique de $(G',{\cal G}',{\bf a})$. La donn\'ee est dite elliptique si $Z(\hat{G}')^{\Gamma_{k},0}=Z(\hat{G})^{\theta,\Gamma_{k},0}$. 
Deux donn\'ees endoscopiques $(G'_{1},{\cal G}'_{1},\tilde{s}_{1})$ et $(G'_{2},{\cal G}'_{2},\tilde{s}_{2})$ sont dites \'equivalentes s'il existe un \'el\'ement $x\in \hat{G}$ tel que $x\tilde{s}_{1}x^{-1}\in\tilde{s}_{2}Z(\hat{G})$ et $x{\cal G}'_{1}x^{-1}={\cal G}'_{2}$. 

 \section{Transformation de la d\'efinition}\label{predonnees}
  Consid\'erons un groupe r\'eductif connexe $\hat{G}$ d\'efini sur ${\mathbb C}$, muni d'une action continue de $\Gamma_{k}$ par automorphismes alg\'ebriques  (continue signifiant qu'elle se quotiente par une action d'un groupe fini $\Gamma_{k'/k}$). On suppose  qu'il existe une paire de Borel \'epingl\'ee $\hat{{\cal E}}=(\hat{B},\hat{T},(E_{\alpha})_{\alpha\in \Delta})$  de $\hat{G}$ qui est conserv\'ee par cette action. On fixe une telle paire de Borel \'epingl\'ee $\hat{{\cal E}}$. On pose encore $^LG=\hat{G}\rtimes W_{k}$. Consid\'erons un automorphisme alg\'ebrique $\theta$ de $\hat{G}$ qui conserve $\hat{{\cal E}}$, est d'ordre fini et commute \`a l'action galoisienne. On a de nouveau des espaces tordus $\hat{G}\theta$ sous $\hat{G}$ et $^LG\theta$ sous $^LG$.  Soit $\tilde{s}$ un \'el\'ement semi-simple de $\hat{G}\theta$. Posons $\hat{G}'=Z_{\hat{G}}(\tilde{s})^0$. Consid\'erons un sous-groupe ${\cal G}'\subset {^LG} $. On suppose ${\cal G}'\cap \hat{G}=\hat{G}'$. On a donc une suite 
$$(1) \qquad 1\to \hat{G}'\to {\cal G}'\to W_{k}\to 1.$$
On suppose que cette suite est exacte et scind\'ee,  au m\^eme sens que dans le paragraphe pr\'ec\'edent.  On suppose que

(2)  pour tout $(g,w)\in {\cal G}'$, il existe $z\in Z(\hat{G})$ de sorte que $\tilde{s}(g,w)=(zg,w)\tilde{s}$. 

A ces conditions, appelons le couple $({\cal G}',\tilde{s})$ une pr\'e-donn\'ee endoscopique de $(^LG,\theta)$. Les notions d'ellipticit\'e d'une pr\'e-donn\'ee ou d'\'equivalence de deux pr\'e-donn\'ees se d\'efinissent comme dans le paragraphe pr\'ec\'edent.

Dans la situation de ce paragraphe pr\'ec\'edent, une donn\'ee endoscopique de $(G,\tilde{G},{\bf a})$ d\'efinit une pr\'e-donn\'ee de $(^LG,\theta)$. Inversement, soit $({\cal G}',\tilde{s})$ une pr\'e-donn\'ee. A isomorphisme pr\`es, il existe un unique groupe r\'eductif connexe $G'$ d\'efini et quasi-d\'eploy\'e sur $k$ de sorte que $\hat{G'}$ en soit le groupe dual. Il existe alors une unique classe ${\bf a}$ de sorte que le triplet $(G', {\cal G}',\tilde{s})$ soit une donn\'ee endoscopique de $(G,\tilde{G},{\bf a})$. Ces correspondances
conservent la notion  d'\'equivalence. Autrement dit, l'ensemble des  pr\'e-donn\'ees endoscopiques de $(^LG,\theta)$ s'identifie \`a la r\'eunion sur les classes ${\bf a}$ des ensembles de donn\'ees endoscopiques de $(G,\tilde{G},{\bf a})$.

 Remarquons que la notion de pr\'e-donn\'ee ne concerne que des objets d\'efinis sur ${\mathbb C}$. Le corps $k$ n'y intervient que via une action continue  de $\Gamma_{k}$.  
 
 \section{Rappels concernant les pr\'e-donn\'ees endoscopiques}\label{rappels}
 On se place dans la situation du paragraphe pr\'ec\'edent.  On note  $w\mapsto w_{G}$ l'action de $W_{k}$ sur $\hat{G}$ qui \'etend  l'action galoisienne. 
 
 On note $\Sigma$ l'ensemble des racines de $\hat{T}$ dans $\hat{G}$. Pour $\alpha\in \Sigma$, on note $\alpha_{res}$ la restriction de $\alpha$ \`a $\hat{T}^{\theta,0}$. On pose $\Sigma_{res}=\{\alpha_{res}; \alpha\in \Sigma\}$. C'est un syst\`eme de racines en g\'en\'eral non r\'eduit. On note $\Sigma_{res,ind}$ le sous-ensemble des racines indivisibles, qui est encore un syst\`eme de racines. L'automorphisme $\theta$ agit sur $\Sigma$ et $\Sigma_{res}$ est en bijection avec l'ensemble des orbites pour cette action. Notons $\Delta_{res}=\{\alpha_{res}; \alpha\in \Delta\}$. Cet ensemble est une base de $\Sigma_{res,ind}$. 
 Pour $\alpha\in \Sigma$, on note $e(\alpha)$ le plus petit $n\in {\mathbb N}_{>0}$ tel que $\theta^n(\alpha)=\alpha$. Pour $\beta\in \Sigma_{res}$ et pour une racine $\alpha\in \Sigma$ telle que $\beta=\alpha_{res}$, l'entier $e(\alpha)$ ne d\'epend que de $\beta$, on le note aussi $e(\beta)$. 
 
 Notons $W$ le groupe de Weyl de $\hat{G}$ relativement \`a $\hat{T}$. L'automorphisme $\theta$  agit sur $W$. Le groupe $W^{\theta}$ est le groupe de Weyl du syst\`eme de racines $\Sigma_{res,ind}$. Plus pr\'ecis\'ement, $W^{\theta}$ s'identifie au groupe de Weyl du groupe $\hat{G}^{\theta,0}$ relativement au sous-tore maximal $\hat{T}^{\theta,0}$. Le syst\`eme de racines de $\hat{G}^{\theta,0}$ est $\Sigma_{res,ind}$, cf. \cite{KS} 1.3.  
 
 Soit $({\cal G}',\tilde{s})$ une pr\'e-donn\'ee endoscopique. A \'equivalence pr\`es, on peut supposer $\tilde{s}=s\theta$, avec $s\in \hat{T}^{\theta,0}$.  Supposons cette condition v\'erifi\'ee. On a not\'e $\hat{G}'=Z_{\hat{G}}(\tilde{s})^0$. Le tore $\hat{T}^{\theta,0}$ en est un sous-tore maximal. Notons $W^{\hat{G}'}$ le groupe de Weyl de $\hat{G}'$ relatif \`a ce tore. Alors $W^{\hat{G}'}\subset W^{\theta}$. 
 Notons $\Sigma(\hat{G}')$ l'ensemble de racines de $\hat{T}^{\theta,0}$ dans $\hat{G}'$.  D'apr\`es \cite{KS} 1.3,  $\Sigma(\hat{G}')$ est le sous-ensemble des racines $\beta\in \Sigma_{res}$ v\'erifiant l'\'egalit\'e
 
 $\beta(s)^{e(\beta)}=1$ si $\beta$ est indivisible;
 
 $\beta(s)^{e(\beta)}=-1$ si $\beta$ est divisible.
 
 Le couple $(\hat{B}\cap \hat{G}',\hat{T}^{\theta,0})$ est une paire de Borel de $\hat{G}'$. Compl\'etons-la en une paire de Borel \'epingl\'ee $\hat{{\cal E}}'$. Le groupe ${\cal G}'$ normalise $\hat{G}'$. Pour tout $w\in W_{k}$, on peut donc fixer $ (g(w),w)\in {\cal G}'$ tel que la restriction \`a $\hat{G}'$ de l'automorphisme  $ad_{g(w)}\circ w_{G}$ conserve $\hat{{\cal E}}'$. On note $w\mapsto w_{G'}$ l'action de $W_{k}$ sur $\hat{G}'$, $w\in W_{k}$ agissant par la restriction de l'automorphisme pr\'ec\'edent. Puisque $ad_{g(w)}\circ  w_{G}$  conserve $\hat{{\cal E}}'$  et que $w_{G}$ conserve $\hat{T}^{\theta,0}$, 
  la conjugaison par $g(w)$ conserve $\hat{T}^{\theta,0}$, donc aussi son commutant $\hat{T}$. Donc $g(w)$ appartient \`a $Norm_{\hat{G}}(\hat{T})$. Notons $u(w)$ son image dans $W$. Son action sur $\hat{T}$ conserve $\hat{T}^{\theta,0}$ et cela implique que $u(w)$ appartient \`a  $W^{\theta}$. D'apr\`es ce que l'on a dit ci-dessus, on peut fixer un \'el\'ement $n(w)\in Norm_{\hat{G}^{\theta,0}}(\hat{T}^{\theta,0})$ dont l'image dans $W^{\theta}$ soit $u(w)$. Il existe alors $t(w)\in \hat{T}$ tel que $g(w)=t(w)n(w)$. On a l'\'egalit\'e ${\cal G}'=\{ht(w)n(w),w); w\in W_{k}, h\in \hat{G}'\}$.  La condition   \ref{predonnees}(2)  \'equivaut \`a
 
 (1) pour tout $w\in W_{k}$, il existe $z\in Z(\hat{G})$ de sorte que $w_{G'}(s)=zs\theta(t(w))t(w)^{-1}$.

  L'existence d'une section continue de la projection ${\cal G}'\mapsto W_{k}$ entra\^{\i}ne que les applications $w\mapsto  n(w),w_{G'}$  peuvent \^etre choisies telles qu'elles se quotientent par un quotient fini de $W_{k}$, puis s'\'etendent en des applications d\'efinies sur $\Gamma_{k}$. On notera celles-ci $\sigma\mapsto n(\sigma),\sigma_{G'}$.

 Introduisons quelques notations.  On note $1-\theta$ l'homomorphisme $t\mapsto t\theta(t)^{-1}$ de $\hat{T}$ dans lui-m\^eme. On note $\hat{G}_{AD}$ le groupe adjoint de $\hat{G}$. Pour un \'el\'ement $g\in \hat{G}$, resp. un sous-groupe $\hat{H}$ de $\hat{G}$, on note $g_{ad}$, resp.$\hat{H}_{ad}$, son image dans $\hat{G}_{AD}$. Signalons que les groupes $\hat{T}_{ad}^{\theta}$ et $\hat{G}_{AD}^{\theta}$ sont connexes.

\section{Les pr\'e-donn\'ees endoscopiques elliptiques sont d'ordre fini}\label{ordrefini}

On dit qu'une pr\'e-donn\'ee endoscopique est $s$-unitaire  si elle est \'equivalente 
  \`a une donn\'ee $({\cal G}',\tilde{s})$ telle que
  $\tilde{s}=s\theta$, o\`u $s$ est un \'el\'ement de  $ \hat{T}^{\theta,0}$ tel que $\beta(s)\in U^1$ pour tout $\beta\in \Sigma_{res}$ (rappelons que $U^1$ est le groupe des nombres complexes de module $1$). On dit qu'une pr\'e-donn\'ee endoscopique est d'ordre fini si elle \'equivalente \`a une  donn\'ee $({\cal G}',\tilde{s})$ telle que
  $\tilde{s}=s\theta$, o\`u $s$ est un \'el\'ement de  $ \hat{T}^{\theta,0}$ 
   dont l'image dans $\hat{T}^{\theta,0}/Z(\hat{G})^{\theta,0}$ est   d'ordre fini. Une pr\'e-donn\'ee d'ordre fini est $s$-unitaire.
  
  {\bf Remarque.} On v\'erifie ais\'ement  qu'une pr\'e-donn\'ee endoscopique  est $s$-unitaire, resp. d'ordre fini, si et seulement si, pour toute pr\'e-donn\'ee $({\cal G}',\tilde{s})$  qui lui est \'equivalente, telle que $\tilde{s}=s\theta$, avec $s\in \hat{T}^{\theta,0}$,  $s$ v\'erifie la condition requise ci-dessus: $\beta(s)\in U^1$ pour toute $\beta\in \Sigma_{res}$, resp.  l'image de $s$ dans $\hat{T}^{\theta,0}/Z(\hat{G})^{\theta,0}$ est   d'ordre fini. 
  
\begin{lem}{Toute pr\'e-donn\'ee endoscopique elliptique est d'ordre fini, a fortiori $s$-unitaire.}\end{lem}

Preuve.  Nous adaptons au cas tordu la preuve de Langlands, cf. \cite{L} p. 705. On consid\`ere une pr\'e-donn\'ee $({\cal G}',\tilde{s})$ avec $\tilde{s}=s\theta$, o\`u $s\in \hat{T}^{\theta,0}$.  Supposons que l'image de $s$ dans $\hat{T}^{\theta,0}/Z(\hat{G})^{\theta,0}$ n'est pas  d'ordre fini. On va prouver que la pr\'e-donn\'ee n'est pas elliptique. On utilise les d\'efinitions et notations du paragraphe \ref{rappels}. 
 Notons $\Phi$ l'ensemble des $\beta\in \Sigma_{res,ind}$ tels que $\beta(s)\in \mu({\mathbb C})$.  D'apr\`es la description donn\'ee au paragraphe  \ref{rappels} de l'ensemble $\Sigma(\hat{G}')$, tout \'el\'ement de cet ensemble est multiple   d'un \'el\'ement de $\Phi$.  L'hypoth\`ese que l'image de $s$ dans $\hat{T}^{\theta,0}/Z(\hat{G})^{\theta,0}$ n'est pas d'ordre fini signifie que  $\Phi$ n'est pas \'egal \`a $\Sigma_{res,ind}$ tout entier. Notons $\pi:\Sigma_{res,ind}\to {\mathbb C}^{\times}/\mu({\mathbb C})$ l'application qui, \`a $\beta\in \Sigma_{res,ind}$ associe l'image de $\beta(s)$ dans ${\mathbb C}^{\times}/\mu({\mathbb C})$. Notons $A$ le groupe engendr\'e par l'image de $\pi$.  Alors $A$ est un sous-groupe de type fini de ${\mathbb C}^{\times}/\mu({\mathbb C})$, non r\'eduit \`a l'\'el\'ement neutre. Puisque ${\mathbb C}^{\times}/\mu({\mathbb C})$ est sans torsion, $A$ est un groupe ab\'elien libre et on peut en fixer une base $\{a_{1},...,a_{n}\}$. On a $n\geq1$. Notons $A^M$ le sous-groupe de $A$ engendr\'e par $\{a_{1},...,a_{n-1}\}$ et $A^P$ le sous-ensemble de $A$ form\'e des \'el\'ements dont la $n$-i\`eme coordonn\'ee dans la base $\{a_{1},...,a_{n}\}$. est positive ou nulle. 
Notons $\Sigma_{res,ind}^M$, resp. $\Sigma_{res,ind}^P$, l'ensemble des $\beta\in \Sigma_{res,ind}$ telles que $\pi(\beta)\in A^M$, resp. $\pi(\beta)\in A^P$. L'ensemble $\Sigma_{res,ind}^P$ est un sous-ensemble parabolique propre de $\Sigma_{res,ind}$ et $\Sigma_{res,ind}^M$ en est le sous-ensemble de Levi. L'ensemble $\Phi$ est inclus dans $\Sigma_{res,ind}^M$. Notons $Z^M_{*}\subset X_{*}(\hat{T}^{\theta,0})$ l'annulateur  de $\Sigma_{res,ind}^M$, c'est-\`a-dire l'ensemble des $x_{*}\in X_{*}(\hat{T}^{\theta,0})$ tels que $\beta(x_{*})=0$ pour tout $\beta\in \Sigma_{res,ind}^M$. Puisque tout \'el\'ement de $\Sigma(\hat{G}')$ est multiple d'un \'el\'ement de $\Phi$, lequel est contenu dans $ \Sigma_{res,ind}^M$, $Z^M_{*}$ est contenu dans $X_{*}(Z(\hat{G}')^0)$. 

L'action galoisienne $\sigma\mapsto \sigma_{G'}$ de $W_{k}$ sur $\hat{T}^{\theta,0}$ conserve $\Sigma_{res,ind}$.  Montrons que

(1) cette action conserve $\Sigma_{res,ind}^M$ et $\Sigma_{res,ind}^P$. 

Soit $\sigma\in \Gamma_{k}$. Il suffit de prouver que $\pi\circ \sigma_{G'}^{-1}=\pi$. Posons $s'=\sigma_{G'}(s)$. Pour $\beta\in \Sigma_{res,ind}$, on a $\sigma_{G'}^{-1}(\beta)(s)=\beta(s')$. En vertu de la d\'efinition de $\pi$, il suffit de prouver que $\beta(s')\in \mu({\mathbb C})\beta(s)$ pour tout $\beta\in \Sigma_{res,ind}$.  La relation \ref{rappels}(1) entra\^{\i}ne que  $s'_{ad}$ appartient \`a $(1-\theta)(\hat{T}_{ad})s_{ad}$. Puisque $s_{ad}$ et $s'_{ad}$ appartiennent tous deux \`a $\hat{T}_{ad}^{\theta}$, on a m\^eme $s'_{ad}\in s_{ad} \left((1-\theta)(\hat{T}_{ad})\cap \hat{T}_{ad}^{\theta}\right)$. Mais ce groupe $(1-\theta)(\hat{T}_{ad})\cap \hat{T}_{ad}^{\theta}$ est fini. Il en r\'esulte que, pour tout $\beta\in \Sigma_{res,ind}$, le quotient $\beta(s')\beta(s)^{-1}$ est une racine de l'unit\'e. C'est la propri\'et\'e voulue qui d\'emontre (1). 

A toute racine $\beta\in \Sigma_{res,ind}$ est associ\'ee une coracine $\check{\beta}\in X_{*}(\hat{T}^{\theta,0})$. Notons $\check{\delta}$ la somme des coracines associ\'ees aux \'el\'ements $\beta\in \Sigma_{res,ind}^P-\Sigma_{res,ind}^M$. C'est un \'el\'ement de $Z^M_{*}$ qui n'appartient pas \`a $X_{*}(Z(\hat{G})^0)$. La propri\'et\'e (1) entra\^{\i}ne que $\check{\delta}$ est fix\'e par l'action $\sigma\mapsto \sigma_{G'}$. Puisque $Z^M_{*}\subset X_{*}(Z(\hat{G'})^0)$, on obtient que $\check{\delta}$ est un \'el\'ement de $X_{*}(Z(\hat{G}')^{\Gamma_{k},0})$ qui n'appartient pas \`a $X_{*}(Z(\hat{G})^{\theta,\Gamma_{k},0})$. Donc la pr\'e-donn\'ee $({\cal G}',\tilde{s})$ n'est pas elliptique. Cela ach\`eve la preuve. $\square$

 \section{Constructions immobili\`eres}\label{constructionsimmobilieres}

Pour ce paragraphe et les suivants jusqu'en \ref{equivalencepresquepartout}, on suppose que $\hat{G}$ est adjoint.  

Notons $e$ l'ordre de $\theta$. Fixons un nombre premier $p>e$ et une puissance $q$ de $p$ telle que $q\equiv 1\,mod\,e{\mathbb Z}$. Fixons un corps local $F$ non-archim\'edien de caract\'eristique nulle dont le corps r\'esiduel est ${\mathbb F}_{q}$. On peut fixer une extension galoisienne totalement ramifi\'ee $E/F$ de sorte que $\Gamma_{E/F}\simeq {\mathbb Z}/e{\mathbb Z}$. On fixe un g\'en\'erateur $\gamma\in \Gamma_{E/F}$. Il existe \`a isomorphisme pr\`es un unique groupe r\'eductif connexe ${\bf G}$ d\'efini et quasi-d\'eploy\'e sur $F$ de sorte que les conditions suivantes soient v\'erifi\'ees. Fixons une paire de Borel \'epingl\'ee $\boldsymbol{{\cal E}}=({\bf B},{\bf T},({\bf E}_{\alpha})_{\alpha\in \Delta})$ de ${\bf G}$ conserv\'ee par l'action galoisienne. Alors il existe des isomorphismes en dualit\'e $X_{*}(\hat{T})\simeq X_{*}({\bf T})$ et $X^{*}(\hat{T})\simeq X^{*}({\bf T})$ qui font se correspondre les coracines, resp. racines, de $\hat{T}$ dans $\hat{G}$ et de ${\bf T}$ dans ${\bf G}$, ainsi que les coracines, resp. racines, simples relativement \`a $\hat{B}$, resp. ${\bf B}$ (c'est pourquoi nous avons not\'e par anticipation $\Delta$ l'ensemble commun de racines simples).  Par ces isomorphismes, l'action de $\theta$ sur $\Delta$ correspond \`a l'action de $\gamma\in \Gamma_{E/F}$. On notera encore $\theta$ l'action alg\'ebrique sur ${\bf T}$ d\'eduite de l'action de $\gamma$, ou encore transport\'ee par les isomorphismes pr\'ec\'edents de l'action de $\theta$ sur $\hat{T}$. Le groupe ${\bf G}$ est adjoint et de m\^eme type que $\hat{G}$. On note  encore $\Sigma$ l'ensemble commun de   racines de ${\bf T}$ dans ${\bf G}$ ou de $\hat{T}$ dans $\hat{G}$. De l'action du  groupe $\Gamma_{k}$ sur $\hat{G}$ se d\'eduit une action sur $\Sigma$ qui pr\'eserve $\Delta$. Elle se transporte donc en une action de $\Gamma_{k}$ sur ${\bf G}$ par automorphismes alg\'ebriques qui conservent $\boldsymbol{{\cal E}}$ et commutent \`a l'action de $\gamma$. On note $\sigma\mapsto \sigma_{G}$ cette action.

Introduisons les immeubles de Bruhat-Tits $Imm_{E}({\bf G})$ et $Imm_{F}({\bf G})$  de ${\bf G}$ sur $E$ et $F$. Au tore ${\bf T}$ est associ\'e un appartement $App_{E}({\bf T})$ dans $Imm_{E}({\bf G})$. A 
la paire de Borel \'epingl\'ee $\boldsymbol{{\cal E}}$  est associ\'ee une alc\^ove  $C\subset App_{E}({\bf T})$.  

Supposons un instant que ${\bf G}$ est simple (sur $E$). 
Notons  $\alpha_{0}$ l'oppos\'ee de la racine de plus grande longueur de $\Sigma$ (relativement \`a la base $\Delta$) et posons $\Delta_{a}=\Delta\cup \{\alpha_{0}\}$. Notons $\bar{C}$ l'adh\'erence de $C$. Alors $\Delta_{a}$ s'identifie \`a l'ensemble des  sommets du simplexe $\bar{C}$. On note $\alpha\mapsto s_{\alpha}$ cette bijection. L'appartement $App_{E}({\bf T})$ s'identifie \`a $X_{*}({\bf T})\otimes_{{\mathbb Z}}{\mathbb R}$, le sommet $s_{\alpha_{0}}$ s'identifiant \`a $0$.  Revenons au cas g\'en\'eral o\`u ${\bf G}$ n'est pas suppos\'e simple. Notons  
$${\bf G}=\prod_{i\in I({\bf G})}{\bf G}_{i}$$
la d\'ecomposition de ${\bf G}$ en composantes simples. Les objets introduits ci-dessus qui ne d\'ependent que de la structure de ${\bf G}$ sur $E$ ont des analogues pour chaque ${\bf G}_{i}$ et on affecte ces analogues d'un indice $i$. On a 
$$App_{E}({\bf T})=\prod_{i\in I({\bf G})}App_{E}({\bf T}_{i})\simeq \prod_{i\in I({\bf G})}X_{*}({\bf T}_{i})\otimes_{{\mathbb Z}}{\mathbb R}\simeq X_{*}({\bf T})\otimes_{{\mathbb Z}}{\mathbb R},$$
et $C=\prod_{i\in I({\bf G})}C_{i}$. On pose $\Delta_{a}=\sqcup_{i\in I({\bf G})}\Delta_{a,i}$.

Notons ${\bf T}^{nr}$ le plus grand sous-tore de ${\bf T}$ qui est d\'eploy\'e sur $F$. On a $X_{*}({\bf T}^{nr})=X_{*}({\bf T})^{\Gamma_{E/F}}=X_{*}({\bf T})^{\theta}$. Le groupe $\Gamma_{E/F}$ agit sur $Imm_{E}({\bf G})$ et $Imm_{F}({\bf G})$ s'identifie \`a l'ensemble des points fixes de cette action (parce que l'extension $E/F$ est mod\'er\'ement ramifi\'ee). En particulier, l'appartement $App_{F}({\bf T}^{nr})$ s'identifie \`a l'ensemble des points fixes de l'action  naturelle de $\theta$ dans $App_{E}({\bf T})$, ou encore \`a $X_{*}({\bf T}^{nr})\otimes_{{\mathbb Z}}{\mathbb R}$. Posons 
$${\bf A}^{nr}=App_{F}({\bf T}^{nr})=X_{*}({\bf T}^{nr})\otimes_{{\mathbb Z}}{\mathbb R}.$$
Notons $C^{nr}$ l'ensemble des points fixes dans $C$ pour l'action de $\theta$. C'est une alc\^ove dans ${\bf A}^{nr}$.  Toute racine  $\beta\in \Sigma_{res}$ d\'efinit une forme lin\'eaire sur $ {\bf A}^{nr}$. Pour $\lambda\in {\mathbb R}$, notons $\beta[\lambda]$ la fonction affine $x\mapsto \beta(x)+\lambda$ sur ${\bf A}^{nr}$. Bruhat et Tits d\'efinissent l'ensemble des racines affines de cet appartement. Pour tout $\beta\in \Sigma_{res}$, il existe un sous-ensemble discret  $\Lambda(\beta)\subset {\mathbb Q}$ de sorte que l'ensemble des racines affines soit \'egal \`a $\{\beta[\lambda]; \beta\in \Sigma_{res},\lambda\in \Lambda(\beta)\}$. Rappelons la description de cet ensemble, cf. \cite{KP} 2.2.7.    
 On a les \'egalit\'es
 
 $\Lambda(\beta)=\frac{1}{e(\beta)}{\mathbb Z}$, si $\beta$ est indivisible;
 
  $\Lambda(\beta)=\frac{1}{e(\beta)}(\frac{1}{2}+{\mathbb Z})$, si $\beta$ est divisible.
  
  Notons $\Theta$ le groupe d'automorphismes de ${\bf G}$ engendr\'e par $\theta$. Ce groupe $\Theta$ agit naturellement sur $I({\bf G})$. On note $I({\bf G})/\Theta$ l'ensemble d'orbites. Supposons un instant qu'il n'y a qu'une orbite. Pour  $\beta\in \Delta_{res}$, on a $0\in \Lambda(\beta)$ et on pose $\beta^{aff}=\beta[0]$. Il existe une unique racine $\beta_{0}\in \Sigma_{res}$ et un unique \'el\'ement $\lambda_{0}\in \Lambda(\beta_{0})$ tels qu'en posant $\beta_{0}^{aff}=\beta_{0}[\lambda_{0}]$ et $\Delta_{a}^{nr}=\Delta_{res}\cup\{\beta_{0}\}$, $C^{nr}$ soit le sous-ensemble des $x\in {\bf A}^{nr}$ tels que $\beta^{aff}(x)>0$ pour tout $\beta\in \Delta_{a}^{nr}$. De plus, il existe  d'uniques coefficients r\'eels $(d_{\beta})_{\beta\in \Delta_{a}^{nr}}$ de sorte que l'on ait l'\'egalit\'e
$$(1) \qquad \sum_{\beta\in \Delta_{a}^{nr}}d_{\beta}\beta^{aff}(x)=1$$
pour tout $x\in {\bf A}^{nr}$. Ces coefficients sont en fait  des entiers strictement positifs. L'adh\'erence $\bar{C}$ est un simplexe dont les sommets sont en bijection avec $\Delta_{a}^{nr}$, on note $\beta\mapsto s_{\beta}$ cette bijection. Revenons au cas g\'en\'eral o\`u l'on ne suppose plus que $I({\bf G})/\Theta$ est r\'eduit \`a un unique \'el\'ement. Pour $j\in I({\bf G})/\Theta$, posons ${\bf G}_{j}=\prod_{i\in j}{\bf G}_{i}$. Ce groupe est d\'efini sur $F$ et les objets introduits ci-dessus conservent un sens pour lui. On les affecte d'indices $j$. On a ${\bf A}^{nr}=\prod_{j\in I({\bf G})/\Theta}{\bf A}^{nr}_{j}$, $C^{nr}=\prod_{j\in I({\bf G})/\Theta}C^{nr}_{j}$ et on pose $\Delta_{a}^{nr}=\sqcup_{j\in I({\bf G})/\Theta}\Delta_{a,j}^{nr}$. 

Les objets intervenant dans les descriptions ci-dessus  sont ind\'ependants des corps $F$ et $E$.  

 Le groupe $W$ agit lin\'eairement sur $App_{E}({\bf T})$ et $W^{\theta}$ agit lin\'eairement sur ${\bf A}^{nr}$. Bruhat et Tits d\'efinissent un groupe de Weyl affine $W^{aff}$ qui agit par transformations affines sur ${\bf A}^{nr}$ et qui est produit semi-direct de $W^{\theta}$ et d'un groupe de translations que nous notons $R$. Ce groupe de translations est un sous-groupe de ${\bf A}^{nr}$. Il contient le groupe $X_{*}({\bf T}^{nr}$) comme sous- groupe  d'indice fini. On note $p_{W^{\theta}}:W^{aff}\to W^{\theta}$ la projection naturelle. Remarquons que $W^{\theta}$ intervient \`a la fois comme sous-groupe et comme quotient de $W^{aff}$. 
 
  En fait, on peut d\'efinir diff\'erents groupes de Weyl affines et celui que nous consid\'erons est le plus gros possible. Il contient  comme sous-groupe d'indice fini un sous-groupe $W^{aff}_{sc}$ qui est engendr\'e par les sym\'etries associ\'ees aux racines affines mais, en g\'en\'eral, il ne lui est pas \'egal. On note $\Omega^{nr}$ le sous-groupe des \'el\'ements de $W^{aff}$ qui conservent l'alc\^ove $C^{nr}$. C'est un groupe  ab\'elien fini.   On a la d\'ecomposition en produit semi-direct $W^{aff}=W^{aff}_{sc}\rtimes \Omega^{nr}$, cf. \cite{HR} lemme 14. Introduisons les diagrammes de Dynkin ${\cal D}^{nr}$ et ${\cal D}^{nr}_{a}$ associ\'es aux ensembles de racines $\Delta_{res}$ et $\Delta_{a}^{nr}$. Notons $Aut({\cal D}^{nr})$ et $Aut({\cal D}_{a}^{nr})$ leurs groupes finis d'automorphismes.  Le groupe $Aut({\cal D}^{nr})$ est \'egal au sous-groupe des \'el\'ements de $Aut({\cal D}^{nr}_{a})$ qui permutent  les points du diagramme associ\'es aux  racine $\beta_{0,j}$ pour $j\in I({\bf G})/\Theta$.  Le groupe  $\Omega^{nr}$ s'identifie \`a un sous-groupe de $Aut({\cal D}_{a}^{nr})$. Il s'av\`ere que $Aut({\cal D}_{a}^{nr})=\Omega^{nr}\rtimes Aut({\cal D}^{nr})$. 

Notons $Out({\bf G})$ l'ensemble des automorphismes alg\'ebriques de ${\bf G}$ qui conservent $\boldsymbol{{\cal E}}$. Notons $Out({\bf G})^{\theta}$ le sous-ensemble de ceux qui commutent \`a $\theta$. Le groupe $Out({\bf G})$ agit sur $App_{E}({\bf T})$ en conservant $C$. Il s'en d\'eduit une action sur le diagramme de Dynkin ${\cal D} $ associ\'e \`a l'ensemble de racines simples $\Delta$. En fait, $Out({\bf G})$ s'identifie au groupe des automorphismes de ${\cal D}$. 
Le groupe  $Out({\bf G})^{\theta}$ permute les groupes ${\bf G}_{j}$ pour $j\in I({\bf G})/\Theta$. Il agit sur ${\bf A}^{nr}$  en conservant $C^{nr}$ et $\Delta^{nr}$.  On en d\'eduit un homomorphisme $Out({\bf G})^{\theta}\to Aut({\cal D}^{nr})$. Cette application est bijective si $\theta=1$.   En tout cas, on obtient un homomorphisme $\Gamma_{k}\to Aut({\cal D}^{nr})\subset Aut({\cal D}^{nr}_{a})$, que l'on note encore $\sigma\mapsto \sigma_{ G}$. 

Soit $x\in {\bf A}^{nr}$. Notons $\Sigma^{aff}(x)$ l'ensemble des racines affines $\beta[\lambda]$ telles que $\beta[\lambda](x)=0$. Notons $W_{sc}(x)$ le sous-groupe de $W_{sc}^{aff}$ engendr\'e par les sym\'etries associ\'ees aux \'el\'ements de $\Sigma^{aff}(x)$. L'action de $W_{sc}(x)$ sur ${\bf A}^{nr}$ fixe $x$. Supposons $x\in \bar{C}^{nr}$. Montrons que

(2)  l'application $v\mapsto v(C^{nr})$ est une bijection de $W_{sc}(x)$ sur l'ensemble des alc\^oves $C^{nr}_{\dag}$ de ${\bf A}^{nr}$ telles que $x\in \bar{C}^{ nr}_{\dag}$.

Preuve. Bruhat et Tits ont associ\'e \`a $x$ un groupe r\'eductif connexe ${\bf G}_{x}$ d\'efini sur ${\mathbb F}_{q}$. Du tore ${\bf T}^{nr}$ est issu un sous-tore maximal ${\bf T}_{x}$ de ${\bf G}_{x}$ qui est d\'efini sur ${\mathbb F}_{q}$ et maximalement d\'eploy\'e. En fait, nos hypoth\`eses impliquent que ${\bf G}_{x}$ est d\'eploy\'e:  l'action galoisienne de $\Gamma_{E/F}$ sur ${\bf G}$ se fait par $\theta$ et devient triviale sur ${\bf T}^{nr}$. En particulier, ${\bf T}_{x}$ est d\'eploy\'e. 
Le groupe $W_{sc}(x)$ s'identifie au  groupe de Weyl de ${\bf G}_{x}$ relativement \`a ${\bf T}_{x}$. Notons ${\cal B}_{x}$ l'ensemble des sous-groupes de Borel de ${\bf G}_{x}$ qui contiennent ${\bf T}_{x}$. L'ensemble des alc\^oves $C_{\dag}^{nr}$ de ${\bf A}^{nr}$ telles que $x\in \bar{C}_{\dag}^{nr}$ correspond bijectivement  \`a ${\cal B}_{x}$. En particulier, $C^{nr}$ correspond \`a un groupe ${\bf B}_{x}\in {\cal B}_{x}$. On sait bien que l'application $v\mapsto v({\bf B}_{x})$ est une bijection du groupe de Weyl de ${\bf G}_{x}$ sur ${\cal B}_{x}$. Cela prouve (2). $\square$

Supposons $x\in \bar{C}^{nr}$. Notons $S(x)$ l'ensemble des $\beta\in \Delta_{a}^{nr}$ telles que $\beta^{aff}(x)=0$ et $S^{aff}(x)=\{\beta^{aff}; \beta\in S(x)\}$. On a

(3) $\Sigma^{aff}(x)$ est l'ensemble des racines affines qui sont combinaisons lin\'eaires \`a coefficients dans ${\mathbb Z}$ d'\'el\'ements de $S^{aff}(x)$.

Preuve. Il est clair que toute telle combinaison lin\'eaire annule $x$. Inversement, on introduit le groupe ${\bf G}_{x}$ et ses sous-groupes ${\bf T}_{x}$ et ${\bf B}_{x}$  de la preuve pr\'ec\'edente. D'apr\`es Bruhat et Tits, l'ensemble $\Sigma^{aff}(x)$ s'identifie \`a celui des racines de ${\bf T}_{x}$ dans ${\bf G}_{x}$ tandis que $S^{aff}(x)$ s'identifie \`a celui des racines simples relativement au Borel ${\bf B}_{x}$, cf. \cite{T} 3.5.1, 3.5.2. Donc tout \'el\'ement de $\Sigma^{aff}(x)$ est combinaison lin\'eaire \`a coefficients dans ${\mathbb Z}$ d'\'el\'ements de $S^{aff}(x)$. $\square$

   \section{L'ensemble $Endo$}\label{lensembleendo}
 
  Consid\'erons l'ensemble $\underline{Endo}$ des couples $(\omega_{\star},x)$ tels que
 
 $\omega_{\star}$ est une application continue de $\Gamma_{k}$ dans $\Omega^{nr}$ telle que l'application $\sigma\mapsto \sigma_{\star}:=\omega_{\star}(\sigma)\sigma_{G}$ soit un homomorphisme (\`a valeurs dans le groupe d'automorphismes de ${\bf A}^{nr}$);
 
 $x$ est un \'el\'ement de $\bar{C}^{nr}$;
 
 on a $\sigma_{\star}(x)=x$ pour tout $\sigma\in \Gamma_{k}$.
 
 Deux \'el\'ements $(\omega_{\star},x)$ et $(\omega_{\star'},x')$ de $\underline{Endo}$ sont dits \'equivalents s'il existe $\omega\in \Omega^{nr}$ de sorte que $\omega(x')=x$ et $\sigma_{\star'}=\omega\sigma_{\star}\omega^{-1}$ pour tout $\sigma\in \Gamma_{k}$. On note $Endo$ l'ensemble des classes d'\'equivalence. 
 
 Le groupe $\Theta\times \Gamma_{k}$ agit naturellement sur $I({\bf G})$ et $\Gamma_{k}$ agit sur $I({\bf G})/\Theta$. Notons $I({\bf G})/(\Theta\times \Gamma_{k})$ l'ensemble commun d'orbites et $i({\bf G},\theta\times \Gamma_{k})$ son nombre d'\'el\'ements. 
  Soit $(\omega_{\star},x)\in \underline{Endo}$. De  l'action $\sigma\mapsto \sigma_{\star}$ de $\Gamma_{k}$ sur ${\bf A}^{nr}$  se d\'eduit une action alg\'ebrique de $\Gamma_{k}$ sur ${\bf T}^{nr}$ (on remplace $\omega_{\star}(\sigma)$ par $p_{W^{\theta}}(\omega_{\star}(\sigma))$). On note cette action $\sigma\mapsto \sigma_{\star,alg}$. Cette action conserve 
  $\Delta_{a}^{nr}$ et $S(x)$. Montrons que
  
  (1) le nombre d'orbites de cette action dans $\Delta_{a}^{nr}-S(x)$ a au moins $i({\bf G},\theta\times \Gamma_{k})$ \'el\'ements.
  
   Notons $\xi:\Delta_{a}^{nr}\to I({\bf G})/\Theta$ l'application qui, \`a $\beta\in \Delta_{a}^{nr}$, associe l'\'el\'ement $j\in I({\bf G})/\Theta$ tel que $\beta\in \Delta_{a,j}^{nr}$.  Munissons $\Delta_{a}^{nr}$ de l'action $\sigma\mapsto \sigma_{\star,alg}$ de $\Gamma_{k}$ et $I({\bf G})/\Theta$ de l'action $\sigma\mapsto \sigma_{G}$. Alors $\xi$ est \'equivariante pour ces actions.   D\'ecomposons $x$ en $x=\prod_{j\in I({\bf G})/\Theta}x_{j}$ o\`u $x_{j}\in \bar{C}^{nr}_{k}$. Pour tout $j$, l'ensemble $\Delta_{a,j}^{nr}-S(x_{j})$ n'est pas vide.  Donc la restriction de $\xi$ \`a $\Delta_{a}^{nr}-S(x)$ est surjective. L'assertion (1) en r\'esulte.

 On dit que $(\omega_{\star},x)$ est elliptique si  le nombre d'orbites de l'action $\sigma\mapsto \sigma_{\star,alg}$ dans $\Delta_{a}^{nr}-S(x)$  est \'egal \`a $i({\bf G},\theta\times \Gamma_{k})$.   On note $\underline{Endo}_{ell}$ l'ensemble des couples elliptiques et $Endo_{ell}$ l'ensemble de leurs classes d'\'equivalence.

 \section{Un lemme pr\'eliminaire}\label{unlemmepreliminaire}

  Puisque $\hat{G}$ est adjoint, $\hat{T}^{\theta}$ est connexe et \'egal \`a $X_{*}(\hat{T})^{\theta}\otimes_{{\mathbb Z}}{\mathbb C}^{\times}$.   Notons $\hat{T}^{\theta}_{u}$ le sous-groupe des \'el\'ements unitaires de $\hat{T}^{\theta}$, c'est-\`a-dire $\hat{T}^{\theta}_{u}=X_{*}(\hat{T})^{\theta}\otimes_{{\mathbb Z}}U^1$.  L'homomorphisme $z\mapsto e^{2\pi iz}$ identifie ${\mathbb R}/{\mathbb Z}$ \`a $U^1$. 
  On en d\'eduit des isomorphismes
$$\hat{T}^{\theta}_{u}\simeq (X_{*}(\hat{T})^{\theta}\otimes_{{\mathbb Z}}{\mathbb R})/ X_{*}(\hat{T})^{\theta}\simeq (X_{*}({\bf T}^{nr})\otimes_{{\mathbb Z}} {\mathbb R})/X_{*}({\bf T}^{nr})={\bf A}^{nr}/X_{*}({\bf T}^{nr}).$$
 L'intersection $(1-\theta)(\hat{T})\cap \hat{T}^{\theta}$ est finie, contenue dans $\hat{T}^{\theta}_{u}$.  On a introduit au paragraphe  \ref{constructionsimmobilieres} un sous-groupe $R$ de ${\bf A}^{nr}$, qui contient $X_{*}({\bf T}^{nr})$. 

\begin{lem}{L'isomorphisme ci-dessus $\hat{T}^{\theta}_{u}\simeq {\bf A}^{nr}/X_{*}({\bf T}^{nr})$ envoie $(1-\theta)(\hat{T})\cap \hat{T}^{\theta}$ sur $R/X_{*}({\bf T}^{nr})$.}\end{lem}

 Preuve du lemme.   Notons $val_{F}$ la valuation usuelle de $F$ et prolongeons-la en une valuation de $\bar{F}$ \`a valeurs dans ${\mathbb Q}\cup\{\infty\}$. On a $val_{F}(E^{\times})=\frac{1}{e}{\mathbb Z}$. Pour $\alpha\in \Delta$, notons $E_{\alpha}$ le sous-corps de $E$ form\'e des points fixes par $\gamma^{e(\alpha)}$. On a $[E_{\alpha}:F]=\frac{1}{e(\alpha)}$ et $val_{F}(E_{\alpha}^{\times})=\frac{1}{e(\alpha)}{\mathbb Z}$. 
 
 De la valuation $val_{F}$ se d\'eduit un homomorphisme
$$val_{{\bf T}}:{\bf T}(\bar{F})=X_{*}({\bf T})\otimes_{{\mathbb Z}}\bar{F}^{\times}\to X_{*}({\bf T})\otimes_{{\mathbb Z}}{\mathbb Q} .$$
Cet homomorphisme est \'equivariant pour les actions de $\Gamma_{F}$, donc l'image de $ {\bf T}(F)={\bf T}(\bar{F})^{\Gamma_{F}}$ est contenue dans ${\bf A}^{nr}$. D'apr\`es Bruhat et Tits, $R$ est \'egal \`a cette image, cf. \cite{T} 1.2. Introduisons la base $(\check{\varpi}_{\alpha})_{\alpha\in \Delta}$ de $X_{*}({\bf T})$ duale de $\Delta$. L'automorphisme $\theta$ agit encore sur cette base. Pour $\alpha\in \Delta$, on pose $\check{\varpi}_{\alpha_{res}}=\sum_{n=1,...,e(\alpha)}\check{\varpi}_{\theta^n(\alpha)}$. Alors $(\check{\varpi}_{\beta})_{\beta\in \Delta_{res}}$ est la base de $X_{*}({\bf T}^{nr})$ duale de $\Delta_{res}$. Puisque ${\bf T}$ est d\'eploy\'e sur $E$, on a ${\bf T}(E)=X_{*}({\bf T})\otimes_{{\mathbb Z}}E^{\times}$. Soit $t=\prod_{\alpha\in \Delta}\check{\varpi}_{\alpha}(t_{\alpha})\in {\bf T}(E)$, avec des  coefficients $t_{\alpha}\in E^{\times}$. On a $t\in {\bf T}(F)$ si et seulement si $\gamma(t_{\alpha})=t_{\theta(\alpha)}$ pour tout $\alpha\in  \Delta$. Cela entra\^{\i}ne $t_{\alpha}\in E_{\alpha}^{\times}$. Inversement, la condition peut \^etre r\'ealis\'ee pour tout 
  $t_{\alpha}\in E_{\alpha}^{\times}$, en supposant que  $t_{\theta^n(\alpha)}=\gamma^n(t_{\alpha})$ pour $n=1,...,e(\alpha)-1$. Il s'en d\'eduit que l'image $val_{{\bf T}}({\bf T}(F))$, c'est-\`a-dire $R$, est \'egale \`a $\oplus_{\beta\in \Delta_{res}}\frac{1}{e_{\beta}}{\mathbb Z}\check{\varpi}_{\beta}$. Soit maintenant $z=\prod_{\alpha\in \Delta}\check{\varpi}_{\alpha}(z_{\alpha})$ un \'el\'ement de $\hat{T}$, avec des $z_{\alpha}\in {\mathbb C}^{\times}$. Il appartient \`a $(1-\theta)(\hat{T})$ si et seulement s'il existe $z'=\prod_{\alpha\in \Delta}\check{\varpi}_{\alpha}(z'_{\alpha})$  tel que $z=(1-\theta)(z')=
 \prod_{\alpha\in \Delta}\check{\varpi}_{\alpha}(z'_{\alpha}z_{\theta^{-1}(\alpha)}^{'-1})$. Cela \'equivaut \`a ce que, pour tout $\alpha\in \Delta$, on ait $\prod_{n=0,...,e(\alpha)-1} z_{\theta^n(\alpha)}=1$. D'autre part, $z$ appartient \`a $\hat{T}^{\theta}$ si et seulement si $z_{\theta^n(\alpha)}=z_{\alpha}$ pour tout $\alpha\in \Delta$ et tout $n=1,...,e(\alpha)-1$. Donc $z\in (1-\theta)(\hat{T})\cap \hat{T}^{\theta}$ si et seulement si, outre cette derni\`ere  condition, on a de plus $z_{\alpha}\in \mu_{e(\alpha)}({\mathbb C})$. Autrement dit, $t=\prod_{\beta\in \Delta_{res}}\check{\varpi}_{\beta}(z_{\beta})$ avec des $z_{\beta}\in \mu_{e(\beta)}({\mathbb C})$. En comparant avec le calcul de $R$, on obtient le lemme. $\square$ 
 
 Il se d\'eduit du lemme un isomorphisme
 $$\iota:\hat{T}_{u}^{\theta}/((1-\theta)(\hat{T})\cap \hat{T}^{\theta})\to {\bf A}^{nr}/R.$$
 Les groupes $W^{\theta}$ et $\Gamma_{k}$ agissent naturellement sur ces deux groupes et $\iota$ entrelace ces actions.  Notons $p_{{\bf A}^{nr}}:{\bf A}^{nr}\to {\bf A}^{nr}/R$ la projection naturelle. Le groupe $W^{aff}$ agit sur ${\bf A}^{nr}$ et $W^{\theta}$ agit sur ${\bf A}^{nr}/R$. On a la relation $p_{{\bf A}^{nr}}(v(x))=(p_{W^{\theta}}(v))(p_{{\bf A}^{nr}}(x))$ pour tous $x\in {\bf A}^{nr}$ et $v\in W^{aff}$. On note $p_{\hat{T}^{\theta}}:\hat{T}_{u}^{\theta}\to \hat{T}_{u}^{\theta}/((1-\theta)(\hat{T})\cap \hat{T}^{\theta})$ la projection naturelle.

 \section{Construction de pr\'e-donn\'ees endoscopiques}\label{constructiondepredonnees}

 Soit $(\omega_{\star},x)\in \underline{Endo}$. Fixons $s\in \hat{T}_{u}^{\theta}$ tel que $\iota\circ p_{\hat{T}^{\theta}}(s)=p_{{\bf A}^{nr}}(x)$ et posons $\tilde{s}=s\theta$.  Pour $\sigma\in \Gamma_{k}$, posons $u(\sigma)=p_{W^{\theta}}(\omega_{\star}(\sigma))$. Relevons l'application $\sigma\mapsto u(\sigma)$ en une application continue $n:\Gamma_{k}\to Norm_{\hat{G}^{\theta}}(\hat{T}^{\theta})$.  Pour $\sigma\in \Gamma_{k}$, l'\'egalit\'e $\sigma_{\star}(x)=x$ et les propri\'et\'es d'entrelacement de l'isomorphisme $\iota$ et de nos diff\'erentes projections impliquent l'\'egalit\'e $u(\sigma)\sigma_{G}(s)\in s(1-\theta)(\hat{T})$. Donc $n(\sigma)\sigma_{G}(s)n(\sigma)^{-1}\in s(1-\theta)(\hat{T})$. Fixons une application continue $t:\Gamma_{k}\to \hat{T}$ de sorte que $n(\sigma)\sigma_{G}(s)n(\sigma)^{-1}=st(\sigma)^{-1}\theta(t(\sigma))$ pour tout $\sigma\in \Gamma_{k}$. Posons $\tilde{s}=s\theta$ et  $g(\sigma)=t(\sigma)n(\sigma)$ pour tout $\sigma\in \Gamma_{k}$. Les applications $\sigma\mapsto n(\sigma),g(\sigma)$ etc... se prolongent contin\^ument au groupe de Weil $W_{k}$.
 Parce que $n(w)\in \hat{G}^{\theta}$, l'\'egalit\'e $n(w)w_{G}(s)n(w)^{-1}=st(w )^{-1}\theta(t(w))$  \'equivaut \`a 
 
 (1) $\tilde{s}(g(w),w)=(g(w),w)\tilde{s}$ pour tout $w\in W_{k}$. 
 
 Posons $\hat{G}'=Z_{\hat{G}}(\tilde{s})^0$ et d\'efinissons ${\cal G}'=\{hg(w),w); w\in W_{k}, h\in \hat{G}'\}$. Cet ensemble ne d\'epend pas des choix de $n(w)$ et $t(w)$. En effet $n(w)$ est uniquement d\'etermin\'e modulo $\hat{T}^{\theta}\subset \hat{G}'$ et $t(w)$ l'est aussi car  l'\'el\'ement $(1-\theta)(t(w))$ est uniquement d\'etermin\'e.
 
   \begin{lem}{Le couple $({\cal G}',\tilde{s})$ est une pr\'e-donn\'ee endoscopique $s$-unitaire.}\end{lem}
 
 Preuve.  On voit que les \'el\'ements $n(1)$ et $t(1)$ appartiennent \`a  $ \hat{T}^{\theta}$, donc \`a $ \hat{G}'$.  D'o\`u  aussi $g(1)\in \hat{G}'$. Cela entra\^{\i}ne ${\cal G}'\cap \hat{G}=\hat{G}'$. 
 D'apr\`es (1), on a la relation cl\'e 
 
 (2) $\tilde{s}(g,w)=(g,w)\tilde{s}$ pour tout $(g,w)\in {\cal G}'$. 
 
 Il reste \`a prouver que ${\cal G}'$ est un groupe et que la projection ${\cal G}'\to W_{k}$ admet une section qui soit un homomorphisme continu.  Le fait que l'application $\sigma\mapsto \sigma_{\star}$ est un homomorphisme implique que, pour $\sigma,\sigma'\in \Gamma_{k}$, on a $n(\sigma)\sigma(n(\sigma'))\in \hat{T}^{\theta}n(\sigma\sigma')$. Il en r\'esulte que, pour $w,w'\in W_{k}$, on a la relation  $g(w)w(g(w'))\in \hat{T}g(ww')$.  Ecrivons $g(w)w(g(w'))=t_{2} (w,w')g(ww')$, avec $t_{2}(w,w')\in \hat{T}$. Parce que $(g(w),w)$, $(g(w'),w')$ et $(g(ww'),ww')$ v\'erifient tous les trois la relation (2), on a $t_{2}(w,w')\in \hat{T}^{\theta}$. Alors $(g(w),w)(g(w'),w')=(t_{2}(w,w') g(ww'),ww')$ appartient \`a ${\cal G}'$. Un argument analogue prouve que $(g(w),w)^{-1}$ appartient \`a ${\cal G}'$.   D'autre part, la relation (2) implique que tout \'el\'ement de ${\cal G}'$ normalise $\hat{G}'$.   De ces  propri\'et\'es r\'esulte que ${\cal G}'$ est un groupe.    L'application $w\mapsto (g(w),w)$ est une section continue de la projection ${\cal G}'\to W_{k}$. Ce n'est pas forc\'ement un homomorphisme. D\'efinissons une action de $\Gamma_{k}$ sur $\hat{T}$ par $\sigma\mapsto  \sigma_{G'}=u(\sigma)\circ \sigma_{G}$. Elle conserve $\hat{T}^{\theta}$. 
 L'application $(w,w')\mapsto t_{2}(w,w')$ d\'efinie ci-dessus est par construction un $2$-cocycle pour cette action, qui provient d'un $2$-cocycle continu d\'efini sur $\Gamma_{k}$. On a vu qu'elle prenait ses valeurs dans $\hat{T}^{\theta}$. Parce que ce groupe  est un tore, ce cocycle est trivial, cf. \cite{L} lemme 4. On peut donc choisir une application continue  $w\mapsto t'(w)$ de $W_{k}$ dans $\hat{T}^{\theta}$ dont le cobord soit $t_{2}$. L'application $w\mapsto (t(w)^{'-1}g(w),w)$ est alors un homomorphisme continu qui scinde la projection  ${\cal G}'\to W_{k}$.  Cela ach\`eve la d\'emonstration. $\square$.
 
 La pr\'e-donn\'ee que l'on a construite d\'epend du choix de $s$. On montrera que sa classe d'\'equivalence n'en d\'epend pas.

\section{Identification du groupe $W^{\hat{G}'}$}\label{unepropriete}
Soit $(\omega_{\star},x)\in \underline{Endo}$. Construisons une donn\'ee $({\cal G}',\tilde{s})$ comme dans le paragraphe pr\'ec\'edent, dont nous utilisons les notations. On a d\'efini le sous-groupe $W_{sc}(x)\subset W^{aff}$ au paragraphe \ref{constructionsimmobilieres}. La projection $p_{W^{\theta}}:W^{aff}\to W^{\theta}$ est injective sur $W_{sc}(x)$. En effet, si $v,v'\in W_{sc}(x)$ v\'erifi\'ent $p_{W^{\theta}}(v)=p_{W^{\theta}}(v')$, il existe $r\in R$ tel que $v'=rv$. Or $v$ et $v'$ fixent $x$ et la seule translation qui fixe un point de $ {\bf A}^{nr}$ est $0$. Donc $r=0$ et $v=v'$. On a aussi d\'efini le sous-groupe $W^{\hat{G}'}\subset W^{\theta}$.

\begin{lem}{Le groupe $W^{\hat{G}'}$ est \'egal \`a l'image de $W_{sc}(x)$ par $p_{W^{\theta}}$.}\end{lem}

Preuve. Le groupe $W_{sc}(x)$ est engendr\'e par les sym\'etries associ\'ees aux racines affines appartenant \`a $\Sigma^{aff}(x)$, cf. paragraphe \ref{constructionsimmobilieres}. Notons $p_{\Sigma_{res}}$ l'application qui, \`a une racine affine $\beta[\lambda]$, associe la racine $\beta\in \Sigma_{res}$. Pour la m\^eme raison que ci-dessus, la restriction \`a  $\Sigma^{aff}(x)$  de cette application est injective . On note $\Sigma(x)$ l'image de cette application. Cet ensemble $\Sigma(x)$ est celui des $\beta\in \Sigma_{res}$ pour lesquelles  $\beta(x)$ appartient \`a $\Lambda(\beta)$. L'image de $W_{sc}(x)$ par  $p_{W^{\theta}}$ est le sous-groupe de $W^{\theta}$ engendr\'e par les sym\'etries associ\'ees aux \'el\'ements de $\Sigma(x)$. D'autre part, $W^{\hat{G}'}$ est le sous-groupe de $W^{\theta}$ engendr\'e par les sym\'etries associ\'ees aux \'el\'ements de $\Sigma(\hat{G}')$. Il suffit donc de prouver

(1) $\Sigma(x)=\Sigma(\hat{G}')$. 

Remarquons que, pour $\beta\in \Sigma_{res}$, l'homomorphisme $\beta^{n_{\beta}}:\hat{T}_{u}^{\theta}\to U^1$ se descend en un homomorphisme encore not\'e $\beta^{e_{\beta}}$ sur $\hat{T}_{u}^{\theta}/((1-\theta)(\hat{T})\cap \hat{T}^{\theta})$. En effet,  soit $\alpha\in \Sigma$ tel que $\alpha_{res}=\beta$. Alors $\beta^{e_{\beta}}$ co\"{\i}ncide sur $\hat{T}^{\theta}$ avec $\prod_{n=1,...,e(\alpha)}\theta^n(\alpha)$. Il est clair que cet homomorphisme est trivial sur $(1-\theta)(\hat{T})$. Cela prouve l'assertion. Via l'isomorphisme $\hat{T}_{u}^{\theta}\simeq {\bf A}^{nr}/X_{*}({\bf T}^{nr})$, $\beta$ correspond \`a la compos\'ee de la racine $\beta:{\bf A}^{nr}\to {\mathbb R}$ et de la projection ${\mathbb R}\to {\mathbb R}/{\mathbb Z}$. Notons $\beta_{mod\,{\mathbb Z}}$ cette compos\'ee. Via la remarque pr\'ec\'edente, l'application $e_{\beta}\beta_{mod\,{\mathbb Z}}$ se descend en un homomorphisme $e_{\beta}\beta_{mod\,{\mathbb Z}}:{\bf A}^{nr}/R\to {\mathbb R}/{\mathbb Z}$. En vertu de la d\'efinition $\iota\circ p_{\hat{T}^{\theta}}(s)=p_{{\bf A}^{nr}}(x)$, on a alors les \'equivalences
$$ \beta(s)^{e(\beta)}=1\iff e_{\beta}\beta(x)\in {\mathbb Z},\,\, \beta(s)^{e(\beta)}=-1\iff e_{\beta}\beta(x)\in \frac{1}{2}+{\mathbb Z}.$$
Il suffit alors de comparer la description de $\Sigma(\hat{G}')$ donn\'ee au paragraphe \ref{rappels} et celle de l'ensemble $\Lambda(\beta)$ donn\'ee au paragraphe \ref{constructionsimmobilieres} pour constater que $\beta\in \Sigma(\hat{G}')$ si et seulement si $\beta(x)$ appartient \`a $\Lambda(\beta)$, autrement dit $\beta\in \Sigma(x)$. Cela prouve (1) et le lemme. $\square$

\section{Equivalences}\label{equivalences}
Soient $({\cal G}'_{1},\tilde{s}_{1})$ et $({\cal G}'_{2},\tilde{s}_{2})$ deux pr\'e-donn\'ees endoscopiques. Notons $Isom(\tilde{s}_{1};\tilde{s}_{2})$ l'ensemble des $g\in \hat{G}$ tels que $g\tilde{s}_{1}g^{-1}=\tilde{s}_{2}$.  Notons $Isom({\cal G}'_{1},\tilde{s}_{1}; {\cal G}'_{2},\tilde{s}_{2})$ l'ensemble des 
$g\in Isom(\tilde{s}_{1};\tilde{s}_{2})$ tels que $g{\cal G}'_{1}g^{-1}={\cal G}'_{2}$. Dans le cas o\`u  les deux pr\'e-donn\'ees sont \'egales, notons-les simplement $({\cal G}',\tilde{s})$. Alors   $Isom(\tilde{s};\tilde{s})=Z_{\hat{G}}(\tilde{s})$ et on note $Aut({\cal G}',\tilde{s})=Isom({\cal G}',\tilde{s};{\cal G}',\tilde{s})$. On a l'inclusion $\hat{G}'\subset Aut({\cal G}',\tilde{s})$. Revenons au cas g\'en\'eral o\`u les pr\'e-donn\'ees ne sont pas suppos\'ees \'egales. Si $Isom(\tilde{s}_{1};\tilde{s}_{2})$, resp $Isom({\cal G}'_{1},\tilde{s}_{1}; {\cal G}'_{2},\tilde{s}_{2})$, est non vide, c'est un torseur \`a gauche  sous le groupe $Z_{\hat{G}}(\tilde{s}_{2})$, resp. $Aut({\cal G}'_{2},\tilde{s}_{2})$. 

 Soient $(\omega_{\star,1},x_{1})$ et $(\omega_{\star,2},x_{2})$ deux \'el\'ements de $\underline{Endo}$. Notons $\Omega^{nr}(x_{1};x_{2})$ l'ensemble des $\omega\in \Omega^{nr}$ tels que $\omega(x_{1})=x_{2}$. Notons $Isom(\omega_{\star,1},x_{1};\omega_{\star,2},x_{2})$ l'ensemble des $\omega\in \Omega^{nr}(x_{1};x_{2})$ tels que $\omega\sigma_{\star,1}\omega^{-1}=\sigma_{\star,2}$ pour tout $\sigma\in \Gamma_{k}$. Dans le cas o\`u $(\omega_{\star,1},x_{1})=(\omega_{\star,2},x_{2})$, notons simplement $(\omega_{\star},x)$ cet \'el\'ement. On pose $\Omega^{nr}(x)=\Omega^{nr}(x;x)$ et $Aut(\omega_{\star},x)=Isom(\omega_{\star},x;\omega_{\star},x)$. Revenons au cas g\'en\'eral o\`u les deux \'el\'ements ne sont pas suppos\'es \'egaux. Si $\Omega^{nr}(x_{1};x_{2})$, resp. $Isom( \omega_{\star,1},x_{1};\omega_{\star,2},x_{2})$, n'est pas vide, c'est un torseur \`a gauche sous le groupe $\Omega^{nr}(x_{2})$, resp $Aut(\omega_{\star,2},x_{2})$.

\begin{prop}{(i) Soient $(\omega_{\star,1},x_{1})$ et $(\omega_{\star,2},x_{2})$ deux \'el\'ements de $\underline{Endo}$ et soient $({\cal G}'_{1},\tilde{s}_{1})$ et $({\cal G}'_{2},\tilde{s}_{2})$ deux pr\'e-donn\'ees endoscopiques. Supposons que, pour $i=1,2$, $({\cal G}'_{i},\tilde{s}_{i})$ soit construite \`a partir de $(\omega_{\star,i},x_{i})$ comme dans le paragraphe \ref{constructiondepredonnees}. Alors il existe une bijection entre les ensembles $\hat{G}'_{2}\backslash Isom(\tilde{s}_{1};\tilde{s}_{2})$ et $\Omega^{nr}(x_{1};x_{2})$, qui se restreint en 
 une bijection entre les ensembles $  \hat{G}'_{2}\backslash Isom({\cal G}'_{1},\tilde{s}_{1};{\cal G}'_{2},\tilde{s}_{2})$ et $Isom(\omega_{\star,1},x_{1};\omega_{\star,2},x_{2})$. 
 
 (ii) Soient $(\omega_{\star},x)\in \underline{Endo}$ et $({\cal G}',\tilde{s})$ une pr\'e-donn\'ee endoscopique  construite \`a partir de $(\omega_{\star},x)$ comme dans le paragraphe  \ref{constructiondepredonnees}. Alors il existe un isomorphisme de groupes de $\hat{G}'\backslash  Z_{\hat{G}}(\tilde{s})$ sur $\Omega^{nr}(x)$,  qui  se restreint en un isomorphisme de $\hat{G}'\backslash Aut({\cal G}',\tilde{s})$ sur $Aut(\omega_{\star},x)$. }\end{prop}
 
 Preuve. On se place dans la situation de (i).  Pour simplifier la notation, on pose $\boldsymbol{{\cal G}}'_{i}=({\cal G}'_{i},\tilde{s}_{i})$ et ${\bf E}_{i}=(\omega_{\star,i},x_{i})$ pour $i=1,2$.  On reprend les notations du paragraphe  \ref{constructiondepredonnees} en y ajoutant des indices $i$.  Notons $W(\tilde{s}_{1};\tilde{s}_{2})$ l'ensemble des $u\in W^{\theta}$ tels que $u(s_{1})\in s_{2}((1-\theta)(\hat{T})\cap \hat{T}^{\theta})$. Notons $W(\boldsymbol{{\cal G}}'_{1};\boldsymbol{{\cal G}}'_{2})$ l'ensemble des $u\in W(\tilde{s}_{1};\tilde{s}_{2})$ tels que $uu_{1}(\sigma)\sigma_{G}(u)^{-1}\in W^{\hat{G}'_{2}}u_{2}(\sigma)$ pour tout $\sigma\in \Gamma_{k}$. Montrons que
 
 (1) le groupe $W^{\hat{G}'_{2}}$ agit \`a gauche sur $W(\tilde{s}_{1};\tilde{s}_{2})$ en conservant  $W(\boldsymbol{{\cal G}}'_{1};\boldsymbol{{\cal G}}'_{2})$; cette actions est libre si $W(\tilde{s}_{1};\tilde{s}_{2})$ n'est pas vide; 
 il existe une bijection entre les ensembles $ W^{\hat{G}'_{2}}\backslash W(\tilde{s}_{1};\tilde{s}_{2})$ et $\hat{G}'_{2}\backslash Isom(\tilde{s}_{1};\tilde{s}_{2})$, qui se restreint en une bijection de $W^{\hat{G}'_{2}}\backslash W(\boldsymbol{{\cal G}}'_{1};\boldsymbol{{\cal G}}'_{2})$ sur $\hat{G}'_{2}\backslash Isom(\boldsymbol{{\cal G}}'_{1};\boldsymbol{{\cal G}}'_{2})$. 
 
 Posons $N=Norm_{\hat{G}}(\hat{T}^{\theta})$. Puisque $\hat{T}^{\theta}$ est un sous-tore maximal de $\hat{G}'_{i}$ pour $i=1,2$, on a l'\'egalit\'e $Isom(\tilde{s}_{1};\tilde{s}_{2})=\hat{G}'_{2}(N\cap Isom(\tilde{s}_{1};\tilde{s}_{2}))$.  Il en r\'esulte une bijection
 $$\hat{G}'_{2}\backslash Isom(\tilde{s}_{1};\tilde{s}_{2})\simeq Norm_{\hat{G}'_{2}}(\hat{T}^{\theta})\backslash (N\cap Isom(\tilde{s}_{1};\tilde{s}_{2})).$$
 Deux \'el\'ements de $N\cap Isom(\tilde{s}_{1};\tilde{s}_{2})$, resp. $Norm_{\hat{G}'_{2}}(\hat{T}^{\theta})$,  ont m\^eme image dans $W^{\theta}$ si et seulement s'ils diff\`erent par multiplication \`a gauche par un \'el\'ement de $\hat{T}^{\theta}$. En notant $\underline{W}(\tilde{s}_{1};\tilde{s}_{2})$ l'image de $N\cap Isom(\tilde{s}_{1};\tilde{s}_{2})$ dans $W^{\theta}$, on en d\'eduit que $W^{\hat{G}'_{2}}$ agit \`a gauche sur $\underline{W}(\tilde{s}_{1};\tilde{s}_{2})$, que cette action est libre si cet ensemble n'est pas vide et que l'on a une bijection
 $$\hat{G}'_{2}\backslash Isom(\tilde{s}_{1};\tilde{s}_{2})\simeq  W^{\hat{G}'_{2}}\backslash \underline{W}(\tilde{s}_{1};\tilde{s}_{2}).$$
 De m\^eme, notons $\underline{W}(\boldsymbol{{\cal G}}'_{1};\boldsymbol{{\cal G}}'_{2})$ l'image dans $W^{\theta}$ de $N\cap Isom({\cal G}'_{1};{\cal G}'_{2})$. On a une bijection
 $$ \hat{G}'_{2}\backslash Isom({\cal G}'_{1},\tilde{s}_{1};{\cal G}'_{2},\tilde{s}_{2})\simeq W^{\hat{G}'_{2}}\backslash \underline{W}(\boldsymbol{{\cal G}}'_{1},\boldsymbol{{\cal G}}'_{2}).$$
 Alors (1) r\'esulte des deux assertions suivantes
 
 (2) on a l'\'egalit\'e $W(\tilde{s}_{1};\tilde{s}_{2})=\underline{W}(\tilde{s}_{1};\tilde{s}_{2})$;
 
 (3) on a l'\'egalit\'e $W(\boldsymbol{{\cal G}}'_{1};\boldsymbol{{\cal G}}'_{2})=\underline{W}(\boldsymbol{{\cal G}}'_{1};\boldsymbol{{\cal G}}'_{2})$. 
 
 Soit $u\in W^{\theta}$, relevons $u$ en un \'el\'ement $n\in Norm_{\hat{G}^{\theta}}(\hat{T}^{\theta})$. On a alors $n\tilde{s}_{1}n^{-1}=u(s_{1})\theta$. 
   L'\'el\'ement $u$ appartient \`a  $\underline{W}(\tilde{s}_{1};\tilde{s}_{2})$ si et seulement s'il existe $t\in \hat{T}$ tel que $tn\in  Isom(\tilde{s}_{1};\tilde{s}_{2})$. Cette derni\`ere relation \'equivaut \`a  $tn\tilde{s}_{1}n^{-1}t^{-1}=\tilde{s}_{2}$, ou encore \`a $u(s_{1})t\theta(t)^{-1}=s_{2}$. Donc $u$ appartient \`a $\underline{W}(\tilde{s}_{1};\tilde{s}_{2})$ si et seulement si $u(s_{1})\in s_{2}(1-\theta)(\hat{T})$. Puisque $s_{2}$ et $u(s_{2})$ appartiennent \`a $\hat{T}^{\theta}$, cela \'equivaut \`a $u(s_{1})\in s_{2}((1-\theta)(\hat{T})\cap \hat{T}^{\theta})$, c'est-\`a-dire $u\in W(\tilde{s}_{1},\tilde{s}_{2})$. Cela prouve (2).  
 
 Gr\^ace \`a (2), les deux ensembles $W(\boldsymbol{{\cal G}}'_{1};\boldsymbol{{\cal G}}'_{2})$ et $\underline{W}(\boldsymbol{{\cal G}}'_{1};\boldsymbol{{\cal G}}'_{2})$ sont en tout cas contenus dans $W(\tilde{s}_{1};\tilde{s}_{2})$. Pour prouver (3), on peut supposer que $W(\tilde{s}_{1};\tilde{s}_{2})$ n'est pas vide. Soit $u\in W(\tilde{s}_{1};\tilde{s}_{2})$. Relevons $u$ en un \'el\'ement $n\in  N\cap Isom(\tilde{s}_{1};\tilde{s}_{2})$. Alors $u$ appartient \`a   $\underline{W}(\boldsymbol{{\cal G}}'_{1};\boldsymbol{{\cal G}}'_{2})$  si et seulement s'il existe $t\in \hat{T}$ tel que, d'une part $tn\in  Isom(\tilde{s}_{1};\tilde{s}_{2})$, d'autre part $tng_{1}(\sigma)\sigma_{G}(tn)^{-1}\in \hat{G}_{2}'g_{2}(\sigma)$ pour tout $\sigma\in \Gamma_{k}$. Puisque $n$ lui-m\^eme appartient \`a $ Isom(\tilde{s}_{1},\tilde{s}_{2})$, la premi\`ere condition \'equivaut \`a $t\in \hat{T}^{\theta}$ et la deuxi\`eme condition est insensible \`a la multiplication par un tel \'el\'ement. Celle-ci est donc: pour tout $\sigma\in \Gamma_{k}$, il existe $h(\sigma)\in \hat{G}'_{2}$ tel que $ng_{1}(\sigma)\sigma_{G}(n)^{-1}=h(\sigma)g_{2}(\sigma)$. Puisque tous ces \'el\'ements, except\'e peut-\^etre $h(\sigma)$, normalisent $\hat{T}^{\theta}$, il en est de m\^eme de $h(\sigma)$. En projetant dans $W^{\theta}$, on obtient l'\'egalit\'e $uu_{1}(\sigma)\sigma_{G}(u)^{-1}\in W^{\hat{G}_{2}'}u_{2}(\sigma)$. Donc $u\in W(\boldsymbol{{\cal G}}'_{1};\boldsymbol{{\cal G}}'_{2})$. Inversement, supposons $u\in W(\boldsymbol{{\cal G}}'_{1};\boldsymbol{{\cal G}}'_{2})$. Alors, pour tout $\sigma\in \Gamma_{k}$,  il existe $u'(\sigma)\in Norm_{\hat{G}_{2}'}(\hat{T}^{\theta})$ et $t(\sigma)\in \hat{T}$ de sorte que $ng_{1}(\sigma)\sigma_{G}(n)^{-1}=t(\sigma)u'(\sigma)g_{2}(\sigma)$. La conjugaison par le membre de gauche envoie $\sigma_{G}(s_{2})\theta$ sur $\tilde{s}_{2}$ et il en est de m\^eme de la conjugaison par $u'(\sigma)g_{2}(\sigma)$. Donc la conjugaison par $t(\sigma)$ fixe $\tilde{s}_{2}$, donc $t(\sigma)\in \hat{T}^{\theta}\subset \hat{G}'_{2}$. Alors $ng_{1}(\sigma)\sigma_{G}(n)^{-1}\in \hat{G}_{2}'g_{2}(\sigma)$ et, comme on l'a dit, $u$ appartient \`a  $\underline{W}(\boldsymbol{{\cal G}}'_{1};\boldsymbol{{\cal G}}'_{2})$. Cela prouve (3), d'o\`u (1).  
 
 Notons $W^{aff}(x_{1};x_{2})$ l'ensemble des $v\in W^{aff}$ tels que $v(x_{1})=x_{2}$ et $W^{aff}({\bf E}_{1};{\bf E}_{2})$ le sous-ensemble des $v\in W^{aff}(x_{1};x_{2})$ tels que, pour tout $\sigma\in \Gamma_{k}$, on ait la relation $v\omega_{\star,1}(\sigma)\sigma_{G}(v)^{-1}\in W_{sc}(x_{2})\omega_{\star,2}(\sigma)$.  La projection $p_{W^{\theta}}$ est injective sur $W^{aff}(x_{1};x_{2})$ car le seul \'el\'ement de $R$ qui fixe un point de $ {\bf A}^{nr}$ est $0$.  On a d\'ej\`a prouv\'e au paragraphe  \ref{unepropriete}  que $p_{W^{\theta}}(W^{\hat{G}'_{2}})=W_{sc}(x_{2})$. Montrons que

(4) $p_{W^{\theta}}(W^{aff}(x_{1};x_{2}))=W(\tilde{s}_{1};\tilde{s}_{2})$, $p_{W^{\theta}}(W^{aff}({\bf E}_{1};{\bf E}_{2}))=W(\boldsymbol{{\cal G}}'_{1};\boldsymbol{{\cal G}}'_{2})$.
 
Puisque $\iota\circ p_{\hat{T}^{\theta}}(s_{i})=p_{{\bf A}^{nr}}(x_{i})$ pour $i=1,2$, un \'el\'ement $u$ appartient \`a $W(\tilde{s}_{1};\tilde{s}_{2})$ si et seulement si $u(x_{1})\in x_{2}+R$. Cela \'equivaut \`a ce que $u$ soit la projection dans $W^{\theta}$ d'un \'el\'ement de $ W^{aff}(x_{1};x_{2})$. D'o\`u la premi\`ere \'egalit\'e de (4). Soit $u\in W(\tilde{s}_{1};\tilde{s}_{2})$, relevons $u$ en son unique ant\'ec\'edant   $v\in W^{aff}(x_{1};x_{2})$. Si $v\in W^{aff}({\bf E}_{1};{\bf E}_{2})$, on a $v\omega_{\star,1}(\sigma)\sigma_{G}(v)^{-1}\in W_{sc}(x_{2})\omega_{\star,2}(\sigma)$ pour tout $\sigma\in \Gamma_{k}$. En projetant dans $W^{\theta}$, on obtient $uu_{1}(\sigma)\sigma_{G}(u)^{-1}\in W^{\hat{G}'_{2}}u_{2}(\sigma)$, donc $u$ appartient \`a $W(\boldsymbol{{\cal G}}'_{1};\boldsymbol{{\cal G}}'_{2})$. R\'eciproquement, supposons que $u\in W(\boldsymbol{{\cal G}}'_{1};\boldsymbol{{\cal G}}'_{2})$. En relevant dans $W^{aff}$ la relation pr\'ec\'edente, on voit qu'il existe $r(\sigma)\in R$ et $v_{x}(\sigma)\in W_{sc}(x_{2})$ de sorte que $v\omega_{\star,1}(\sigma)\sigma_{G}(v)^{-1}=r(\sigma)v_{x}(\sigma)\omega_{\star,2}(\sigma)$. Le membre de gauche envoie $\sigma_{G}(x_{2})$ sur $x_{2}$ et il en est de m\^eme de $v_{x}(\sigma)\omega_{\star,2}(\sigma)$. Donc $r(\sigma)$ fixe $x_{2}$, d'o\`u $r(\sigma)=0$. Alors $v$ appartient \`a $W^{aff}({\bf E}_{1};{\bf E}_{2})$, ce qui d\'emontre (4).

Le groupe $W_{sc}(x_{2})$ agit \`a gauche sur $W^{aff}(x_{1};x_{2})$. Cette action pr\'eserve $W^{aff}({\bf E}_{1};{\bf E}_{2})$: l'action $\sigma_{\star,2}$ pr\'eservant $x_{2}$, elle conserve $W_{sc}(x_{2})$. Prouvons que

(5) $W_{sc}(x_{2})\backslash W^{aff}(x_{1};x_{2})\simeq \Omega^{nr}(x_{1};x_{2})$ et $W_{sc}(x_{2})\backslash W^{aff}({\bf E}_{1};{\bf E}_{2})\simeq Isom({\bf E}_{1},{\bf E}_{2})$. 

D'apr\`es \ref{constructionsimmobilieres}(2), on a $\Omega^{nr}\cap W_{sc}(x_{2})=\{1\}$. On voit que 
 $\Omega^{nr}(x_{1};x_{2})=\Omega^{nr}\cap W^{aff}(x_{1};x_{2})$ et $ Isom({\bf E}_{1};{\bf E}_{2})=\Omega^{nr}\cap W^{aff}({\bf E}_{1};{\bf E}_{2})$. Il suffit alors de prouver que $W^{aff}(x_{1};x_{2})\subset W_{sc}(x_{2})\Omega^{nr}$. Soit $v\in W^{aff}(x_{1};x_{2})$. Alors $v$ envoie $C^{nr}$ sur une alc\^ove dont l'adh\'erence contient $x_{2}$. D'apr\`es \ref{constructionsimmobilieres}(2), il existe donc $v_{x}\in W_{sc}(x_{2})$ tel que $v_{x}v$ conserve $C^{nr}$, c'est-\`a-dire $v_{x}v\in \Omega^{nr}$. Cela prouve (5). 
 
 Le (i) de l'\'enonc\'e proposition r\'esulte de (1), (4) et (5). 
 
 Le (ii) est le cas particulier de (i) o\`u 
 $(\omega_{\star,1},x_{1})=(\omega_{\star,2},x_{2})$  et $({\cal G}'_{1},\tilde{s}_{1})=({\cal G}'_{2},\tilde{s}_{2})$. On affirme en plus que les bijections sont des homomorphismes de groupes. Mais c'est clair sur leurs constructions. $\square$

 \section{La classification des pr\'e-donn\'ees endoscopiques  $s$-unitaires pour $\hat{G}$ adjoint}\label{laclassification}
 Notons ${\cal E}ndo_{u}$ l'ensemble des classes d'\'equivalence de pr\'e-donn\'ees endoscopiques $s$-unitaires et notons ${\cal E}ndo_{ell}$ l'ensemble des classes de pr\'e-donn\'ees endoscopiques elliptiques. On a ${\cal E}ndo_{ell}\subset {\cal E}ndo_{u}$ d'apr\`es le lemme \ref{ordrefini}. Dans le paragraphe \ref{constructiondepredonnees}, on a associ\'e une pr\'e-donn\'ee endoscopique $s$-unitaire  \`a tout \'el\'ement de $\underline{Endo}$. Il r\'esulte du (i) de la proposition \ref{equivalences} que cette construction se descend en une application injective $\delta:Endo\to {\cal E}ndo_{u}$. Par abus de notations, on notera aussi $\delta$ l'application compos\'ee $\underline{Endo}\to Endo \stackrel{\delta}{\to }{\cal E}ndo_{u}$.

 \begin{thm}{(i) L'application   $\delta:Endo\to {\cal E}ndo_{u}$ est bijective.
      
 (ii)  L'image de $Endo_{ell}$ par $\delta$ est \'egale \`a ${\cal E}ndo_{ell}$. } \end{thm}
 
 Preuve.   
  Prouvons la surjectivit\'e de $\delta$.   Consid\'erons une pr\'e-donn\'ee endoscopique $s$-unitaire $({\cal G}',\tilde{s})$. On suppose $\tilde{s}=s\theta$, avec $s\in \hat{T}^{\theta}$ et on utilise les constructions et notations du paragraphe \ref{rappels}.  Fixons $y\in {\bf A}^{nr}$ tel que $p_{{\bf A}^{nr}}(y)=\iota\circ p_{\hat{T}^{\theta}}(s)$. Choisissons un \'el\'ement $v\in W^{aff}$ tel que $v(y)\in \bar{C}^{nr}$. Posons $x=v(y)$. Soit $\sigma\in \Gamma_{k}$. La relation \ref{rappels}(2) entra\^{\i}ne que $u(\sigma)\sigma_{G}(y)\in y+R$.  Il existe donc un unique  \'el\'ement $r(\sigma)\in R$  tel que $r(\sigma)u(\sigma)\sigma_{G}(y)=y$. Posons $v(\sigma)=vr(\sigma)u(\sigma)\sigma_{G}(v)^{-1}$. Alors $v(\sigma)\sigma_{G}(x)=x$. L'action de $v(\sigma)\sigma_{G}$ sur ${\bf A}^{nr}$ envoie $C^{nr}$ sur une alc\^ove dont l'adh\'erence contient $x$. D'apr\`es \ref{constructionsimmobilieres}(2), il existe un unique $v_{x}(\sigma)\in W_{sc}(x)$ tel que $v_{x}(\sigma)v(\sigma)\sigma_{G}$ conserve $C^{nr}$. Posons $\omega_{\star}(\sigma)=v_{x}(\sigma)v(\sigma)$. Puisque $v_{x}(\sigma)$ fixe $x$, on a encore $\omega_{\star}(\sigma)\sigma_{G}(x)=x$. Puisque $\omega_{\star}(\sigma)\sigma_{G}$ conserve $C^{nr}$ et que $\sigma_{G}$ conserve aussi cette alc\^ove, $\omega_{\star}(\sigma)$ la conserve aussi, donc $\omega_{\star}(\sigma)\in \Omega^{nr}$. 
 Montrons que
 
 (1) le couple $(\omega_{\star},x)$ appartient \`a $ \underline{Endo}$.

 Les termes $v$ et $y$ \'etant fix\'es,  il y a un unique choix pour les autres  termes. Puisque l'application $\sigma\mapsto u(\sigma)$ est continue, l'application  $\sigma\mapsto \omega_{\star}(\sigma)$ l'est aussi.  Posons $\sigma_{\star}=\omega_{\star}(\sigma)\sigma_{G}$. Il reste \`a prouver que l'application $\sigma\mapsto \sigma_{\star}$ est un homomorphisme.  Soient $\sigma,\sigma'\in \Gamma_{k}$. Par construction, on a $u(\sigma)\sigma_{G}u(\sigma')\sigma'_{G}=u(\sigma\sigma')(\sigma\sigma')_{G}$. Il en r\'esulte qu'il existe $r_{1}(\sigma,\sigma')\in R$ de sorte que  $r(\sigma)u(\sigma)\sigma_{G}r(\sigma')u(\sigma')\sigma'_{G}=r_{1}(\sigma,\sigma')r(\sigma\sigma')u(\sigma\sigma')(\sigma\sigma')_{G}$. Mais   $r(\sigma)u(\sigma)\sigma_{G}r(\sigma')u(\sigma')\sigma'_{G}$ et $r(\sigma\sigma')u(\sigma\sigma')(\sigma\sigma')_{G}$ fixent $y$ et le seul \'el\'ement de $R$ qui fixe un point de l'appartement est $0$. Donc $r_{1}(\sigma,\sigma')=0$ et $r(\sigma)u(\sigma)\sigma_{G}r(\sigma')u(\sigma')\sigma'_{G}=r(\sigma\sigma')u(\sigma\sigma')(\sigma\sigma')_{G}$. Il en r\'esulte que $v(\sigma)\sigma_{G}v(\sigma')\sigma'_{G}=v(\sigma\sigma')(\sigma\sigma')_{G}$. On en d\'eduit que 
 $$(2) \qquad \sigma_{\star}\sigma_{\star}'=v_{x}(\sigma,\sigma')(\sigma\sigma')_{\star},$$
  o\`u $v_{x}(\sigma,\sigma')=v_{x}(\sigma)\sigma_{\star}(v_{x}(\sigma'))v_{x}(\sigma\sigma')^{-1}$. Parce que l'action $\sigma_{\star}$ fixe $x$, elle conserve le groupe $W_{sc}(x)$ canoniquement attach\'e \`a $x$. Il en r\'esulte que $v_{x}(\sigma,\sigma')$ appartient \`a $W_{sc}(x)$. D'autre part, l'\'egalit\'e (2) et les propri\'et\'es de  l'action $\sigma\mapsto \sigma_{\star}$ entra\^{\i}nent que $v_{x}(\sigma,\sigma')$ conserve $C^{nr}$. Or $1$ est l'unique \'el\'ement de $W_{sc}(x)$ qui conserve cette alc\^ove. Donc $v_{x}(\sigma,\sigma')=1$ et la relation (2) montre que l'application $\sigma\mapsto \sigma_{\star}$ est un homomorphisme. Cela prouve (1).
  
  Montrons que
  
  (3)  la classe d'\'equivalence de $({\cal G}',\tilde{s})$ est  \'egale  \`a  $\delta (\omega_{\star},x)$.

  Construisons une pr\'e-donn\'ee endoscopique associ\'ee \`a $(\omega_{\star},x)$. Puisque la notation $({\cal G}',\tilde{s})$ est d\'ej\`a utilis\'ee, on note cette pr\'e-donn\'ee $({\cal G}'_{\star},\tilde{s}_{\star})$ et on affecte d'un indice $\star$ les objets utilis\'es dans sa construction. Reprenons les notations de la construction ci-dessus  de $(\omega_{\star},x)$ et posons $u=p_{W^{\theta}}(v)$, $u_{x}(\sigma)=p_{W^{\theta}}(v_{x}(\sigma))$.  Par d\'efinition, $s_{\star}$ est un \'el\'ement de $\hat{T}^{\theta}_{u}$ tel que $\iota\circ p_{\hat{T}^{\theta}}(s_{\star})=p_{{\bf A}^{nr}}(x)$. On voit que l'on peut choisir $s_{\star}= u(s)$. Soit $\sigma\in \Gamma_{k}$. 
  Par d\'efinition, $u_{\star}(\sigma)$ est l'image de $\omega_{\star}(\sigma)$ par $p_{W^{\theta}}$. On voit que $u_{\star}(\sigma)= v_{x}(\sigma)uu(\sigma)\sigma_{G}(u)^{-1}$. Relevons $v_{x}(\sigma)$ et $u$ en des \'el\'ements $m(\sigma)$ et $n$ appartenant \`a $Norm_{\hat{G}^{\theta}}(\hat{T}^{\theta})$. Par d\'efinition, $n_{\star}(\sigma)$ est un \'el\'ement de ce groupe qui rel\`eve $u_{\star}(\sigma)$. On peut choisir $n_{\star}(\sigma)=m(\sigma)nn(\sigma)\sigma_{G}(n)^{-1}$. Posons $h_{\star}(\sigma)=nt(\sigma)n^{-1}m(\sigma)^{-1}t_{\star}(\sigma)^{-1}$. Par un calcul simple, l'\'egalit\'e pr\'ec\'edente \'equivaut \`a 
  $$(4) \qquad h_{\star}(\sigma)g_{\star}(\sigma)=ng(\sigma)\sigma_{G}(n)^{-1}.$$
   Remarquons que, puisque $n\in \hat{G}^{\theta}$, l'\'egalit\'e $s_{\star}=u(s)$ \'equivaut \`a $n\tilde{s}_{\star} n^{-1}=\tilde{s}$. Les actions $ad_{g(\sigma)}\circ \sigma_{G}$, resp. $ad_{g_{\star}(\sigma)}\circ \sigma_{G}$, fixent $\tilde{s}$, resp. $\tilde{s}_{\star}$. Il r\'esulte de (4) que $ad_{h_{\star}(\sigma)}$ fixe $\tilde{s}_{\star}$. Autrement dit, $h_{\star}(\sigma)\in Z_{\hat{G}}(\tilde{s}_{\star})$. 
On se rappelle que $v_{x}(\sigma)$ appartient \`a $W_{sc}(x)$. En appliquant le lemme \ref{unepropriete}, on  a $u_{x}(\sigma)\in W^{\hat{G}'_{\star}}$. On peut donc relever $u_{x}(\sigma)$ en un \'el\'ement $m_{\star}(\sigma)\in Norm_{\hat{G}'_{\star}}(\hat{T}^{\theta})$. Alors $m(\sigma)\in \hat{T}m_{\star}(\sigma)$ et on en d\'eduit $h_{\star}(\sigma)\in \hat{T}m_{\star}(\sigma)^{-1}$. Soit $\tau(\sigma)\in \hat{T}$ tel que $h_{\star}(\sigma)=\tau(\sigma)m_{\star}(\sigma)^{-1}$.   Puisque $h_{\star}(\sigma)$ et $m_{\star}(\sigma)$ appartiennent \`a $Z_{\hat{G}}(\tilde{s}_{\star})$, il en est de m\^eme de $\tau(\sigma)$. Alors $\tau(\sigma)\in \hat{T}\cap Z_{\hat{G}}(\tilde{s}_{\star})=\hat{T}^{\theta}$. Mais alors $h_{\star}(\sigma)$ appartient \`a $\hat{G}'_{\star}$. L'\'egalit\'e (4), vraie pour tout $\sigma$, implique que ${\cal G}'_{\star}=n{\cal G}'n^{-1}$. Donc $n$ d\'efinit une \'equivalence entre les donn\'ees $({\cal G}',\tilde{s})$ et $({\cal G}'_{\star}, \tilde{s}_{\star})$. Cela prouve (3).

L'assertion (3) entra\^{\i}ne la surjectivit\'e de $\delta$. Cela ach\`eve la preuve du (i) de l'\'enonc\'e. 
     
   Soit $(\omega_{\star},x)\in \underline{Endo}$, construisons une pr\'e-donn\'ee  $({\cal G}',\tilde{s})$ qui lui est associ\'ee. On utilise les notations du paragrephe \ref{constructiondepredonnees}. On identifie $X_{*}(\hat{T}^{\theta})$ \`a $X_{*}({\bf T}^{nr})$.  On dispose de l'action lin\'eaire $\sigma\mapsto \sigma_{\star,alg}$ de $\Gamma_{k}$ sur $X_{*}({\bf T}^{nr})$, cf. paragraphe \ref{constructionsimmobilieres}.  On dispose d'autre part de l'action $\sigma\mapsto \sigma_{G'}$ de $W_{k}$ sur   $X_{*}(\hat{T}^{\theta}) =X_{*}({\bf T}^{nr})$.  Comparons-les. Il faut prendre garde au fait que les termes $g(\sigma)$, $n(\sigma)$ etc... du paragraphe     \ref{constructiondepredonnees} ne sont pas ceux ainsi not\'es associ\'es \`a $({\cal G}',\tilde{s})$  au paragraphe  \ref{rappels}. En effet, l'action de $ad_{g(\sigma)}\circ \sigma_{G}$ pr\'eserve $\hat{G}'$ et $\hat{T}^{\theta}$ mais pas forc\'ement $\hat{B}\cap \hat{G}'$. Mais, pour tout $\sigma\in \Gamma_{k}$, on peut choisir un \'el\'ement $h(\sigma)\in Norm_{\hat{G}'}(\hat{T}^{\theta})$ tel que $ad_{h(\sigma)g(\sigma)}\circ \sigma_{G}$ pr\'eserve ce groupe de Borel. En notant $u_{h}(\sigma)$ la projection de $h(\sigma)$ dans $W^{\hat{G}'}$, il r\'esulte alors des d\'efinitions que l'on a l'\'egalit\'e
   $\sigma_{G'}=u_{h}(\sigma)\sigma_{\star,alg}$. Cela entra\^{\i}ne que les deux actions $\sigma\mapsto \sigma_{G'}$ et $\sigma\mapsto \sigma_{\star,alg}$ co\"{\i}ncident sur $X_{*}(Z(\hat{G}')^0)$. Munissons $X_{*}({\bf T}^{nr})$ de la deuxi\`eme action. Que $({\cal G}',\tilde{s})$ soit elliptique \'equivaut \`a l'\'egalit\'e $X_{*}(Z(\hat{G}')^0)^{\Gamma_{k}}=0$, ou encore \`a  $X^{*}(Z(\hat{G}')^0)^{\Gamma_{k}}=0$.     Le groupe $X_{*}(Z(\hat{G}')^{0})$ est l'annulateur de $\Sigma(\hat{G}')$ dans $X_{*}({\bf T}^{nr})$. En vertu de \ref{unepropriete}(1) et de \ref{constructionsimmobilieres}(3), c'est aussi l'annulateur de $S(x)$. Pour tout ${\mathbb Z}$-module $Y$, notons $Y_{{\mathbb Q}}=Y\otimes_{{\mathbb Z}}{\mathbb Q}$. Pour tout ensemble fini $A$, notons ${\mathbb Q}[A]$ le ${\mathbb Q}$-espace vectoriel de base $A$. Si $A$ est un sous-ensemble lin\'eairement ind\'ependant de $Y$, on a ${\mathbb Q}[A]\subset Y_{{\mathbb Q}}$.  Le groupe  $ X^*_{{\mathbb Q}}(Z(\hat{G}')^0)$ s'identifie \`a   $X^*_{{\mathbb Q}}({\bf T}^{nr})/{\mathbb Q}[S(x)]$.  Pour tout $j\in I({\bf G})/\Theta$, on a une surjection naturelle
    $${\mathbb Q}[\Delta_{a,j}^{nr}]\to X^*_{{\mathbb Q}}({\bf T}_{j}^{nr}).$$
   Son  noyau est la droite port\'ee par l'\'el\'ement  $l_{j}=\sum_{\beta\in \Delta_{a,j}^{nr}}d(\beta)\beta$, cf. \ref{constructionsimmobilieres}(1).  En sommant les applications ci-dessus, on obtient la suite exacte
   $$0\to \oplus_{j\in I({\bf G})/\Theta}{\mathbb Q}[l_{j}]\to {\mathbb Q}[\Delta_{a}^{nr}] \to X^*_{{\mathbb Q}}({\bf T}^{nr})\to 0.$$
   L'action $\sigma\mapsto \sigma_{\star,alg}$  permute les \'el\'ements de $\Delta_{a}^{nr}$ donc d\'efinit une action sur ${\mathbb Q}[\Delta_{a}^{nr}]$. Les homomorphismes de la suite  ci-dessus sont \'equivariants pour les actions galoisiennes.   On d\'eduit de la suite ci-dessus une suite
   $$ 0\to \oplus_{j\in I({\bf G})/\Theta}{\mathbb Q}[l_{j}]\to {\mathbb Q}[\Delta_{a}^{nr}]/{\mathbb Q}[S(x)] \to X^*_{{\mathbb Q}}({\bf T}^{nr})/{\mathbb Q}[S(x)]\to 0.$$
   Cette suite reste exacte. En effet, le seul point \`a v\'erifier est que le premier homomorphisme est injectif. D\'ecomposons $x$ en $x=\prod_{j\in I({\bf G})/\Theta}x_{j}$. On a $S(x)=\sqcup_{j\in I({\bf G})/\Theta}S(x_{j})$ et il suffit de v\'erifier que, pour tout $j$, $l_{j}$ n'appartient pas \`a ${\mathbb Q}[S(x_{j})]$. C'est clair car $\Delta_{a,j}^{nr}-S(x_{j})$ n'est pas vide.  
   
 Que $({\cal G}',\tilde{s})$ soit elliptique \'equivaut \`a l'\'egalit\'e   $(X^*_{{\mathbb Q}}({\bf T}^{nr})/{\mathbb Q}[S(x)])^{\Gamma_{k}}=0$. D'apr\`es la suite ci-dessus, cela \'equivaut \`a l'\'egalit\'e
 $$dim((\oplus_{j\in I({\bf G})/\Theta}{\mathbb Q}[l_{j}])^{\Gamma_{k}})=dim(({\mathbb Q}[\Delta_{a}^{nr}]/{\mathbb Q}[S(x)] )^{\Gamma_{k}}).$$
 L'action de $\Gamma_{k}$ permute les \'el\'ements $l_{j}$ donc la premi\`ere dimension vaut $i({\bf G},\theta\times \Gamma_{k})$. 
 On a l'\'egalit\'e ${\mathbb Q}[\Delta_{a}^{nr}]/{\mathbb Q}[S(x)]={\mathbb Q}[\Delta_{a}^{nr}-S(x)]$ et l'action galoisienne permute les \'el\'ements de $\Delta_{a}^{nr}-S(x)$. En notant $i(x)$ le nombre d'orbites de cette action dans  $\Delta_{a}^{nr}-S(x)$, on a donc $dim(({\mathbb Q}[\Delta_{a}^{nr}]/{\mathbb Q}[S(x)] )^{\Gamma_{k}})=i(x)$. Donc $({\cal G}',\tilde{s})$ est elliptique si et seulement si $i({\bf G},\theta\times \Gamma_{k})=i(x)$, autrement dit si et seulement si   $(\omega_{\star},x)$ est elliptique. Cela prouve le (ii) de l'\'enonc\'e. $\square$.

   \section{Equivalence presque partout}\label{equivalencepresquepartout}
  
   On suppose que $k$ est un corps de nombres.       Une pr\'e-donn\'ee endoscopique $({\cal G}',\tilde{s})$ d\'efinie sur $k$ se localise pour toute place $v\in V_{k}$ en une pr\'e-donn\'ee endoscopique $({\cal G}'_{v},\tilde{s})$ d\'efinie sur $k_{v}$. Notons plus pr\'ecis\'ement $\underline{Endo}_{k}$ et $\underline{Endo}_{k_{v}}$ les ensembles not\'es pr\'ec\'edemment $\underline{Endo}$ relatifs au corps $k$, resp. $k_{v}$. Alors 
    un \'el\'ement $(\omega_{\star},x)\in \underline{Endo}_{k}$ se localise en un \'el\'ement $(\omega_{\star,v},x)$ de $\underline{Endo}_{k_{v}}$. 
    
    Rappelons que l'on suppose que $\hat{G}$ est adjoint. 
    
    \begin{thm}{Soit $({\cal G}'_{1},\tilde{s}_{1})$ et $({\cal G}'_{2},\tilde{s}_{2})$ deux pr\'e-donn\'ees endoscopiques $s$-unitaires d\'efinies sur $k$. Supposons que, pour presque tout $v\in V_{k}$, les pr\'e-donn\'ees $({\cal G}'_{1,v},\tilde{s}_{1})$ et $({\cal G}'_{2,v},\tilde{s}_{2})$ soient \'equivalentes. Alors $({\cal G}'_{1},\tilde{s}_{1})$ et $({\cal G}'_{2},\tilde{s}_{2})$  sont \'equivalentes.} \end{thm}
    
    Preuve.  Gr\^ace au (i) du th\'eor\`eme \ref{laclassification}, l'\'enonc\'e \'equivaut \`a l'assertion suivante:
    
    (1) soient $(\omega_{\star,1},x_{1})$ et $(\omega_{\star,2},x_{2})$ deux \'el\'ements de $\underline{Endo}_{k}$; supposons que, pour presque tout $v\in V_{k}$, $(\omega_{\star,1,v},x_{1})$ et $(\omega_{\star,2,v},x_{2})$ soient \'equivalents; alors $(\omega_{\star,1},x_{1})$ et $(\omega_{\star,2},x_{2})$ sont \'equivalents. 
    
    On va d'abord proc\'eder \`a trois r\'eductions. On a introduit l'ensemble d'orbites $I({\bf G})/(\Theta\times \Gamma_{k})$ de l'action de   $\Gamma_{k}$ dans $I({\bf G})/\Theta$. Pour $J\in I({\bf G})/(\Theta\times \Gamma_{k})$, posons ${\bf G}_{J}=\prod_{j\in J}{\bf G}_{j}$, c'est-\`a-dire que $j$ parcourt les \'el\'ements de $I({\bf G})/\Theta$ dans l'orbite  $J$. Ce groupe est d\'efini sur $F$ et est stable par l'action de $\Gamma_{k}$. On a ${\bf G}=\prod_{J\in  I({\bf G})/(\Theta\times \Gamma_{k})}{\bf G}_{J}$ et tous nos objets se d\'ecomposent  conform\'ement.  Pr\'ecisons nos notations en notant $\underline{Endo}_{k}({\bf G})$ l'ensemble not\'e jusqu'ici $\underline{Endo}_{k}$.  On a alors $\underline{Endo}_{k}({\bf G})=\prod_{J\in  I({\bf G})/(\Theta\times \Gamma_{k})}\underline{Endo}_{k}({\bf G}_{J})$ et on voit qu'il suffit de d\'emontrer l'assertion analogue \`a (1) pour chaque groupe ${\bf G}_{J}$. Cela nous ram\`ene au cas o\`u l'action de $\Gamma_{k}$ dans $I({\bf G})/\Theta$ n'a qu'une seule orbite, ce que l'on suppose d\'esormais.
    
    Fixons un \'el\'ement $j_{0}\in I({\bf G})/\Theta$ et introduisons l'extension $k_{0}$ de $k$ telle que $\Gamma_{k_{0}}$ soit le fixateur de $j_{0}$ dans $\Gamma_{k}$. L'action de $\Gamma_{k_{0}}$ conserve ${\bf G}_{j_{0}}$ et on peut d\'efinir les ensembles $\underline{Endo}_{k_{0}}({\bf G}_{j_{0}})$ et $Endo_{k_{0}}({\bf G}_{j_{0}})$. On a les \'egalit\'es $\bar{C}^{nr}=\prod_{j\in I({\bf G})/\Theta}\bar{C}^{nr}_{j}$ et $\Omega^{nr}=\prod_{j\in I({\bf G})/\Theta}\Omega^{nr}_{j}$. Soit $(\omega_{\star},x)\in \underline{Endo}_{k}({\bf G})$. Posons $x=\prod_{j\in I({\bf G})/\Theta}x_{j}$ et, pour tout $\sigma\in \Gamma_{k}$, $\omega_{\star}(\sigma)=\prod_{j\in I({\bf G})/\Theta}\omega_{\star,j}(\sigma)$. Notons $\omega_{\star,j_{0},k_{0}}$ la restriction de $\omega_{\star,j_{0}}$ \`a $\Gamma_{k_{0}}$. Alors le couple $(\omega_{\star,j_{0},k_{0}},x_{j_{0}})$ appartient \`a $\underline{Endo}_{k_{0}}({\bf G}_{j_{0}})$. On a ainsi d\'efini une application $\underline{Endo}_{k}({\bf G})\to \underline{Endo}_{k_{0}}({\bf G}_{j_{0}})$. Montrons que
    
    (2) cette application se quotiente en une bijection $Endo_{k}({\bf G})\to Endo_{k_{0}}({\bf G}_{0})$. 
    
    Il s'agit d'un cas d'isomorphisme de Shapiro. Nous nous contenterons de donner la construction de l'application en sens oppos\'e $\underline{Endo}_{k_{0}}({\bf G}_{0})\to \underline{Endo}_{k}({\bf G})$, en laissant le lecteur prouver que les deux applications se descendent en des isomorphismes r\'eciproques entre  $Endo_{k}({\bf G})$ et $Endo_{k_{0}}({\bf G}_{0})$. Soit $(\omega_{\star,j_{0},k_{0}},x_{j_{0}})\in \underline{Endo}_{k_{0}}({\bf G}_{j_{0}})$. 
Puisque $I({\bf G})/\Theta$ forme une unique orbite sous l'action de $\Gamma_{k}$, on peut fixer pour tout $j\in I({\bf G})/\Theta$ un \'el\'ement $\sigma_{j}\in \Gamma_{k}$ tel que $\sigma_{j}(j_{0})=j$. Alors $\sigma_{j, G}$ transporte ${\bf G}_{j_{0}}$ et tous les objets qui lui sont associ\'es sur ${\bf G}_{j}$ et tous les m\^emes objets qui lui sont associ\'es. 
 L'ensemble $\{\sigma_{j}; j\in I({\bf G})/\Theta\}$ est un ensemble de repr\'esentants du quotient $\Gamma_{k}/\Gamma_{k_{0}}$. Pour tout $j$, posons $x_{j}=\sigma_{j, G}(x_{j_{0}})$. Posons $x=\prod_{j\in I({\bf G})/\Theta}x_{j}$. Soient $\sigma\in \Gamma_{k}$ et $j\in I({\bf G})/\Theta$. Il existe d'uniques $j'\in I({\bf G})/\Theta$ et $\sigma'\in \Gamma_{k_{0}}$ tels que $\sigma\sigma_{j'}=\sigma_{j}\sigma'$. Posons $\omega_{\star,j}(\sigma)=\sigma_{j,G}(\omega_{\star,j_{0},k_{0}}(\sigma'))$.  C'est un \'el\'ement de $\Omega^{nr}_{j}$. Posons $\omega_{\star}(\sigma)=\prod_{j\in I({\bf G})/\Theta}\omega_{\star,j}(\sigma)$. On v\'erifie que le couple $(\omega_{\star},x)$ appartient \`a $\underline{Endo}_{k}({\bf G})$. Cela d\'efinit l'application cherch\'ee $\underline{Endo}_{k_{0}}({\bf G}_{0})\to \underline{Endo}_{k}({\bf G})$. Comme on l'a dit, on laisse le lecteur achever la preuve de (2). 
    
 Soit  $v_{0}\in V_{k_{0}}$ et notons $v$ la place de $k$ au-dessous de  $v_{0}$. On a encore une application  $\underline{Endo}_{k_{v}}({\bf G})\to \underline{Endo}_{k_{0,v_{0}}}({\bf G}_{j_{0}})$ qui s'inscrit dans le diagramme commutatif
   $$\begin{array}{ccc}\underline{Endo}_{k}({\bf G})&&\to \underline{Endo}_{k_{0}}({\bf G}_{j_{0}})\\ \downarrow&&\downarrow\\ \underline{Endo}_{k_{v}}({\bf G})&\to& \underline{Endo}_{k_{0,v_{0}}}({\bf G}_{j_{0}})\\ \end{array}$$
   L'application $\underline{Endo}_{k_{v}}({\bf G})\to \underline{Endo}_{k_{0,v_{0}}}({\bf G}_{j_{0}})$ ne v\'erifie pas forc\'ement l'analogue de (2) car $I({\bf G})/\Theta$ n'est pas forc\'ement une unique orbite sous $\Gamma_{k_{v}}$. Il est toutefois clair que cette application se quotiente en une application $Endo_{k_{v}}({\bf G})\to Endo_{k_{0,v_{0}}}({\bf G}_{j_{0}})$. Dans la situation de (1), ces consid\'erations impliquent que les  images de $(\omega_{\star,1},x_{1})$ et $(\omega_{\star,2},x_{2})$  par l'application $\underline{Endo}_{k}({\bf G})\to \underline{Endo}_{k_{0}}({\bf G}_{j_{0}})$ sont encore localement \'equivalentes presque partout. Si nous supposons d\'emontr\'ee l'assertion (1) pour l'ensemble $\underline{Endo}_{k_{0}}({\bf G}_{j_{0}})$, ces images sont \'equivalentes et  l'assertion (2) nous dit  que les donn\'ees initiales $(\omega_{\star,1},x_{1})$ et $(\omega_{\star,2},x_{2})$ le sont aussi. On est ramen\'e \`a d\'emontrer (1)  dans le cas o\`u $I({\bf G})/\Theta$ est r\'eduit \`a un unique \'el\'ement, ce que nous supposons d\'esormais. 
   
   Fixons $i_{0}\in I({\bf G})$. Soit $n$ le plus petit entier strictement positif tel que $\theta^n(i_{0})=i_{0}$. C'est un diviseur de l'ordre $e$ de $\theta$. Posons $e_{0}=\frac{e}{n}$. L'automorphisme $\theta^n$ conserve ${\bf G}_{i_{0}}$ et on note $\theta_{i_{0}}$ sa restriction \`a ce groupe. L'automorphisme $\theta_{i_{0}}$ est d'ordre $e_{0}$.  Notons $F_{0}$ le sous-corps de $E$ fix\'e par $\gamma^n$. L'action de $\Gamma_{F_{0}}$ conserve ${\bf G}_{i_{0}}$ et $\gamma^n$  agit sur $X_{*}({\bf T}_{i_{0}})$  comme $\theta_{i_{0}}$. L'appartement ${\bf A}^{nr}$ s'identifie \`a l'appartement similaire relatif au groupe ${\bf G}_{i_{0}}$ muni de son action de $\Gamma_{F_{0}}$. Puisque $I({\bf G})/\Theta$ est  r\'eduit \`a un unique \'el\'ement, on a $I({\bf G})=\{\theta^m(i_{0}); m=0,...,n-1\}$. Puisque 
   l'action de $\Gamma_{k}$ commute \`a $\theta$, il existe un homomorphisme $a:\Gamma_{k}\to {\mathbb Z}/n{\mathbb Z}$ tel que $\sigma_{G}(i_{0})=\theta^{a(\sigma)}(i_{0})$ pour tout $\sigma\in \Gamma_{k}$. L'automorphisme $\sigma_{G}\theta^{-a(\sigma)}$ conserve ${\bf G}_{i_{0}}$. On note $\sigma_{G_{i_{0}}}$ sa restriction \`a ${\bf G}_{i_{0}}$. Alors l'application $\sigma\mapsto \sigma_{G_{i_{0}}}$ munit ${\bf G}_{i_{0}}$ d'une action de $\Gamma_{k}$ qui commute \`a $\theta_{i_{0}}$. Alors les ensembles $\underline{Endo}_{k}({\bf G}_{i_{0}})$ et $Endo_{k}({\bf G}_{i_{0}})$ sont bien d\'efinis et il est imm\'ediat qu'ils sont \'egaux \`a $\underline{Endo}_{k}({\bf G})$ et $Endo_{k}({\bf G})$. Remarquons qu'ici, contrairement \`a la deuxi\`eme r\'eduction ci-dessus, ces \'egalit\'es persistent quand on remplace $k$ par un localis\'e $k_{v}$ puisque c'est le corps $F$ qui a chang\'e, pas le corps $k$. Il nous suffit donc de d\'emontrer l'assertion (1) pour le groupe ${\bf G}_{i_{0}}$. Cela nous ram\`ene au cas o\`u ${\bf G}$ est absolument
 simple, ce que l'on suppose d\'esormais.

       Dans le cas o\`u $\theta=1$, l'assertion (1) r\'esulte de \cite{LW} propositions 1.6 et 2.2.  Nous ne traitons  ici que le cas  $\theta\not=1$.   Le groupe $\Gamma_{k}$ agit sur ${\bf G}$ par des automorphismes qui pr\'eservent $\boldsymbol{{\cal E}}$ et commutent \`a $\theta$. Parce que $\theta\not=1$ un   tel automorphisme est une puissance de $\theta$. En cons\'equence, l'action de $\Gamma_{k}$ est triviale sur tout objet fix\'e par $\theta$, en particulier     
       l'action $\sigma\mapsto \sigma_{G}$ est triviale sur ${\bf A}^{nr}$. Il en r\'esulte que l'ensemble $Isom(\omega_{\star,1},x_{1};\omega_{\star,2},x_{2})$ d\'efini en \ref{equivalences} est celui des $\omega\in \Omega^{nr}(x_{1};x_{2})$ tels que $\omega\omega_{\star,1}(\sigma)\omega^{-1}=\omega_{\star,2}(\sigma)$ pour tout $\sigma\in \Gamma_{k}$. Mais $\Omega^{nr}$ est commutatif.   L'\'equivalence de   $(\omega_{\star,1},x_{1})$ et $(\omega_{\star,2},x_{2})$  \'equivaut donc \`a la r\'eunion des conditions suivantes
  
  (3) $\Omega^{nr}(x_{1},x_{2})\not=\emptyset$;
  
  (4) $\omega_{\star,1}(\sigma)=\omega_{\star,2}(\sigma)$ pour tout $\sigma\in \Gamma_{k}$.
  
  De m\^eme, pour $v\in V_{k}$, l'\'equivalence de   $(\omega_{\star,1,v},x_{1})$ et $(\omega_{\star,2,v},x_{2})$  \'equivaut  \`a la r\'eunion de (3) et de 
  
  (4)(v) $\omega_{\star,1,v}(\sigma)=\omega_{\star,2,v}(\sigma)$ pour tout $\sigma\in \Gamma_{k_{v}}$.
  
  Puisque $\sigma\mapsto \omega_{\star,i}(\sigma)$ est un homomorphisme continu pour $i=1,2$, le th\'eor\`eme de Cebotarev implique que (4) est v\'erifi\'e si et seulement si (4)(v) l'est pour presque tout $v$. Cela prouve (1) et le th\'eor\`eme. $\square$
    
  \section{Le cas presque simple}\label{presquesimple}
  Revenons au cas g\'en\'eral o\`u le corps  $k$ est soit un corps de nombres, soit un corps local de caract\'eristique nulle.  Depuis le paragraphe \ref{constructionsimmobilieres}, on a suppos\'e que $\hat{G}$ \'etait adjoint. Modifions ces hypoth\`eses en supposant que $\hat{G}$ est presque simple, c'est-\`a-dire que $\hat{G}$ est semi-simple et  que $\hat{G}_{AD}$ est simple.   L'automorphisme $\theta$   et l'action de $\Gamma_{k}$ sur $\hat{G}$ se descendent en un automorphisme de $\hat{G}_{AD}$ et une action de $\Gamma_{k}$ sur ce groupe, que l'on note encore $\theta$ et $\sigma\mapsto \sigma_{G}$. 
  Soit $({\cal G}',\tilde{s})$ une pr\'e-donn\'ee endoscopique pour $(\hat{G},\theta)$.  En \'ecrivant $\tilde{s}=s\theta$, on pose $\tilde{s}_{ad}=s_{ad}\theta$. On pose ${\cal G}'_{ad}=\{(g_{ad},w); (g,w)\in {\cal G}'\}$. Alors $({\cal G}'_{ad},\tilde{s}_{ad})$ est une pr\'e-donn\'ee endoscopique pour $ (\hat{G}_{AD},\theta)$.  Il est clair que deux pr\'e-donn\'ees pour $(\hat{G},\theta)$ qui sont \'equivalentes se descendent ainsi en des pr\'e-donn\'ees pour $(\hat{G}_{AD},\theta)$ qui sont \'equivalentes. 
  
  {\bf Remarque.} Si $\theta=1$, toute pr\'e-donn\'ee endoscopique pour $(\hat{G}_{AD},\theta)$ provient par ce proc\'ed\'e d'une pr\'e-donn\'ee endoscopique pour $(\hat{G},\theta)$. Ce n'est pas vrai si $\theta\not=1$. On donnera un contre-exemple au paragraphe \ref{exemples}.  
  
  \begin{thm}{Supposons que $k$ est un corps de nombres.  Si $\hat{G}$ est presque simple,  le th\'eor\`eme \ref{equivalencepresquepartout} est v\'erifi\'e.}\end{thm}
  
  Preuve.  Puisqu'on a d\'ej\`a d\'emontr\'e cette assertion dans \cite{LW} quand $\theta=1$, on suppose $\theta\not=1$. L'examen de tous les syst\`emes de racines possibles montre que, sous cette hypoth\`ese $\theta\not=1$, le groupe $\Omega^{nr}$ a au plus deux \'el\'ements.
  
  Soient  $({\cal G}'_{1},\tilde{s}_{1})$ et $({\cal G}'_{2},\tilde{s}_{2})$ deux pr\'e-donn\'ees endoscopiques $s$-unitaires pour $(\hat{G},\theta)$ d\'efinies sur $k$. On suppose que  $({\cal G}'_{1,v},\tilde{s}_{1})$ et $({\cal G}'_{2,v},\tilde{s}_{2})$ sont \'equivalentes pour presque toute place $v\in V_{k}$. Alors les pr\'e-donn\'ees $({\cal G}'_{1,ad, v},\tilde{s}_{1,ad})$ et $({\cal G}'_{2,ad,v},\tilde{s}_{2,ad})$ sont aussi \'equivalentes pour presque tout $v$. En appliquant le th\'eor\`eme \ref{equivalencepresquepartout}, les donn\'ees $({\cal G}'_{1,ad},\tilde{s}_{1,ad})$ et $({\cal G}'_{2,ad},\tilde{s}_{2,ad})$  sont \'equivalentes. Fixons $g\in \hat{G}$ tel que $g_{ad}$  conjugue  $({\cal G}'_{2,ad},\tilde{s}_{2,ad})$ en $({\cal G}'_{1,ad},\tilde{s}_{1,ad})$. On ne perd rien en rempla\c{c}ant $({\cal G}'_{2},\tilde{s}_{2})$ par sa conjugu\'ee par $g$.  On peut donc supposer
  
  (1) $({\cal G}'_{2,ad},\tilde{s}_{2,ad})=({\cal G}'_{1,ad},\tilde{s}_{1,ad})$. 
  
    On peut aussi supposer que  $\tilde{s}_{1}=s_{1}\theta$, avec $s_{1}\in \hat{T}^{\theta,0}$. Introduisons des objets comme au paragraphe  \ref{rappels} pour  nos deux donn\'ees pr\'e-endoscopiques, que l'on affecte d'indices $1$ et $2$: $g_{1}(\sigma)$, etc... L'\'egalit\'e (1) entra\^{\i}ne que $g_{2,ad}(\sigma)\in \hat{T}_{ad}^{\theta}g_{1,ad}(\sigma)$. Cela entra\^{\i}ne que $g_{2}(\sigma)\in Z(\hat{G})\hat{T}^{\theta,0}g_{1}(\sigma)$. Puisqu'il est loisible de multiplier $g_{2}(\sigma)$ \`a gauche par un \'el\'ement de $\hat{T}^{\theta,0}$, on peut supposer qu'il existe  $z(\sigma)\in Z(\hat{G})$ tel que $g_{2}(\sigma)=z(\sigma)g_{1}(\sigma)$. Ce terme $z(\sigma)$ est bien d\'etermin\'e modulo $Z(\hat{G})\cap \hat{T}^{\theta,0}$. Posons $\mathfrak{Z}= Z(\hat{G})/(Z(\hat{G})\cap \hat{T}^{\theta,0})$ et  notons $\zeta(\sigma)$  l'image de $z(\sigma)$ dans $\mathfrak{Z}$. On obtient une application uniquement d\'etermin\'ee $\sigma\mapsto \zeta(\sigma)$ de $\Gamma_{k}$ dans $ \mathfrak{Z}$. Remarquons que $\mathfrak{Z}$ est muni d'une unique action galoisienne: les actions $\sigma\mapsto \sigma_{G}$ et $\sigma\mapsto \sigma_{G'_{1}}=\sigma_{G'_{2}}$ y co\"{\i}ncident. Puisque ${\cal G}'_{1}$ et ${\cal G}'_{2}$ sont des groupes, on doit avoir $g_{i}(\sigma)\sigma_{G}(g_{i}(\sigma'))\in \hat{T}^{\theta,0}g_{i}(\sigma\sigma')$ pour $i=1,2$ et tous $\sigma,\sigma'\in \Gamma_{k}$. Donc l'application $\zeta$ est un cocycle. On note $\underline{\zeta}$ son image dans $H^1(\Gamma_{k},  \mathfrak{Z})$. Nous ignorons si   cette classe de  cocycle ne d\'epend que de la classe d'\'equivalence de $({\cal G}'_{2},\tilde{s}_{2})$. 
  On v\'erifie toutefois que
    
    (2) soit $z\in Z(\hat{G})$; si on remplace $({\cal G}'_{2},\tilde{s}_{2})$ par $(z{\cal G}'_{2}z^{-1},\tilde{s}_{2})$, la classe de  cocycle $\underline{\zeta}$ ne change pas. 
    
   Posons $Out({\cal G}'_{1,ad},\tilde{s}_{1,ad})=Aut({\cal G}'_{1,ad},\tilde{s}_{1,ad})/\hat{G}'_{1,ad}$.  Fixons un sous-ensemble fini $\Xi\subset \hat{G}$ tel que la famille $\{\xi_{ad}; \xi\in \Xi\}$ soit un ensemble de repr\'esentants dans  $Aut({\cal G}'_{1,ad},\tilde{s}_{1,ad})$ de ce groupe $Out({\cal G}'_{1,ad},\tilde{s}_{1,ad})$.  On a d\'etermin\'e le groupe $Out({\cal G}'_{1,ad},\tilde{s}_{1,ad})$ au paragraphe  \ref{equivalences}. Fixons un \'el\'ement  $(\omega_{\star},x)\in \underline{Endo}$ tel que    la classe de $({\cal G}'_{1,ad},\tilde{s}_{1,ad})$ soit  $\delta(\omega_{\star},x)$. Alors $Out({\cal G}'_{1,ad},\tilde{s}_{1,ad})\simeq Aut(\omega_{\star},x)$.  Pour la m\^eme raison que dans le paragraphe pr\'ec\'edent, l'hypoth\`ese $\theta\not=1$ entra\^{\i}ne l'\'egalit\'e $Aut(\omega_{\star},x)=\Omega^{nr}(x)$. 
  Soit $v\in V_{k}$. Le m\^eme calcul s'applique \`a la pr\'e-donn\'ee $({\cal G}'_{1,v,ad},\tilde{s}_{1,ad})$ : on a $Out({\cal G}'_{1,v,ad},\tilde{s}_{1,ad})\simeq  \Omega^{nr}(x)$, d'o\`u $Out({\cal G}'_{1,v,ad},\tilde{s}_{1,ad})\simeq Out({\cal G}'_{1,ad},\tilde{s}_{1,ad})$. Mais $Aut({\cal G}'_{1,ad},\tilde{s}_{1})$ s'injecte naturellement dans $Aut({\cal G}'_{1,v,ad},\tilde{s}_{1,ad})$. L'isomorphisme ci-dessus implique donc $Aut({\cal G}'_{1,ad},\tilde{s}_{1,ad}) =Aut({\cal G}'_{1,v,ad},\tilde{s}_{1,ad})$. En cons\'equence, la famille $\{\xi_{ad};\xi\in \Xi\}$ est encore un ensemble de repr\'esentants de $Aut({\cal G}'_{1,v,ad},\tilde{s}_{1,ad})/\hat{G}'_{1,ad}$.  
    
    Pour $\xi\in \Xi$, posons ${\cal G}'_{1,\xi}=\xi{\cal G}'_{1}\xi^{-1}$. La pr\'e-donn\'ee $({\cal G}'_{1,\xi},\tilde{s}_{1})$ est \'equivalente \`a $({\cal G}'_{1},\tilde{s}_{1})$ (car $\xi\tilde{s}_{1}\xi^{-1}\in Z(\hat{G})\tilde{s}_{1}$). Elle ne lui est pas forc\'ement \'egale. Mais elle v\'erifie la m\^eme condition que (1), c'est-\`a-dire $({\cal G}'_{1,\xi,ad},\tilde{s}_{1,ad})=({\cal G}'_{1,ad},\tilde{s}_{1,ad})$. On peut lui appliquer les constructions pr\'ec\'edentes, d'o\`u une classe de cocycles $\underline{\zeta}_{\xi}$. On va prouver
    
    (3) il existe $\xi\in \Xi$ tel que $\underline{\zeta}=\underline{\zeta}_{\xi}$. 
    
    Montrons d'abord que cette assertion implique le th\'eor\`eme. On fixe $\xi\in \Xi$ tel que $\underline{\zeta}=\underline{\zeta}_{\xi}$. Cela signifie que l'on peut fixer $z\in Z(\hat{G})$ de sorte que, pour tout $\sigma\in \Gamma_{k}$, on ait $z(\sigma)\in z\sigma(z)^{-1}z_{\xi}(\sigma)$, o\`u $z_{\xi}(\sigma)$ est l'analogue de $z(\sigma)$ pour la pr\'e-donn\'ee $({\cal G}'_{1,\xi},\tilde{s}_{1})$. Alors ${\cal G}'_{2}= z{\cal G}'_{1,\xi}z^{-1}$ et $z$ d\'efinit une \'equivalence de $({\cal G}'_{2},\tilde{s}_{2})$ sur $({\cal G}'_{1,\xi},\tilde{s}_{1})$. Puisque cette derni\`ere pr\'e-donn\'ee est \'equivalente \`a $({\cal G}'_{1},\tilde{s}_{1})$, on obtient bien l'\'enonc\'e. 
    
    Fixons un sous-ensemble $V'_{k}\subset V_{k}$ de compl\'ementaire fini, tel que, pour tout $v\in V'_{k}$, $({\cal G}'_{1,v},\tilde{s}_{1})$ et $({\cal G}'_{2,v},\tilde{s}_{2})$ soient \'equivalentes. Pour prouver (3), on commence par \'etablir l'assertion locale
    
    (4) pour toute place $v\in V'_{k}$, il existe $\xi_{v}\in \Xi$ tel que $\underline{\zeta}_{v}=\underline{\zeta}_{\xi_{v},v}$.

     Pour $v\in V'_{k}$, soit $g_{v}\in \hat{G}$ tel que $g_{v}{\cal G}'_{1,v}g^{-1}={\cal G}'_{2,v}$ et $g_{v}\tilde{s}_{1}g_{v}^{-1}\in Z(\hat{G})\tilde{s}_{2}$. Gr\^ace \`a (1), on a $g_{v,ad}\in Aut({\cal G}'_{1,ad, v},\tilde{s}_{1,ad})$. Il existe donc $\xi_{v}\in \Xi$ tel que $g_{v,ad}\in \xi_{v,ad}\hat{G}'_{1,ad}$, d'o\`u $g_{v}\in Z(\hat{G})\xi_{v}\hat{G}'_{1}$. Il est loisible de multiplier $u_{v}$ \`a droite par un \'el\'ement de $\hat{G}'_{1}$. On peut donc supposer $g_{v}\in Z(\hat{G})\xi_{v}$. En appliquant (2) (qui vaut aussi bien sur le corps de base $k_{v}$), on obtient que $\underline{\zeta}_{v}=\underline{\zeta}_{\xi_{v},v}$.  Cela d\'emontre (4).

    On peut supposer que $1$ appartient \`a $\Xi$.  Les d\'efinitions entra\^{\i}nent que $\underline{\zeta}_{1}=1$, c'est-\`a-dire que $\underline{\zeta}_{1}$ est la classe du cocycle trivial. Pour la m\^eme raison que dans le paragraphe pr\'ec\'edent, l'hypoth\`ese $\theta\not=1$ entra\^{\i}ne que l'automorphisme $\sigma_{G}$ de $\hat{G}$ est une puissance de $\theta$ pour tout $\sigma\in \Gamma_{k}$.    Puisque $\theta$ est d'ordre au plus $3$. Il existe donc une extension cyclique $k_{0}$ de $k$ de degr\'e au plus $3$ de sorte que l'action de $\Gamma_{k_{0}}$ soit triviale.  Montrons que
    
    (5) pour tout $\xi\in \Xi$, l'image de $\underline{\zeta}_{\xi}$ dans $H^1(\Gamma_{k_{0}},\mathfrak{Z})$ est triviale.
    
    On peut supposer $\xi\not=1$. Parce que $Aut(\omega_{\star},x)=\Omega^{nr}(x)$, la proposition \ref{equivalences} entra\^{\i}ne l'\'egalit\'e $Z_{\hat{G}_{AD}}(\tilde{s}_{1,ad})=Aut({\cal G}'_{1,ad},\tilde{s}_{1,ad})$. Le quotient $\hat{G}_{1,ad}'\backslash Z_{\hat{G}_{AD}}(\tilde{s}_{1,ad})\simeq Out({\cal G}'_{1,ad},\tilde{s}_{1,ad})$ est isomorphe \`a un sous-groupe de $\Omega^{nr}$, lequel a au plus deux \'el\'ements. Il en r\'esulte les \'egalit\'es
    $$Z_{\hat{G}_{AD}}(\tilde{s}_{1,ad})=Aut({\cal G}'_{1,ad},\tilde{s}_{1,ad})=\hat{G}'_{1,ad}\sqcup \hat{G}'_{1,ad}\xi_{ad}.$$
    Puisque $\xi_{ad}$ fixe $\tilde{s}_{1,ad}$, $\xi$ normalise $Z_{\hat{G}}(\tilde{s}_{1})$ et aussi sa composante neutre $\hat{G}'_{1}$.    Puisque $s_{1}\in  \hat{T}^{\theta}$, $s_{1}$ est fix\'e par $\Gamma_{k}$. Il en r\'esulte que $g_{1}(\sigma)_{ad}\in Z_{\hat{G}_{AD}}(\tilde{s}_{1,ad})$. Donc $g_{1}(\sigma)\in Z(\hat{G})\hat{G}_{1}'\sqcup Z(\hat{G})\hat{G}_{1}'\xi$. Alors $\xi g_{1}(\sigma)\xi^{-1}\in \hat{G}'g_{1}(\sigma)$. Mais, pour $\sigma\in \Gamma_{k_{0}}$, cela \'equivaut \`a $\xi g_{1}(\sigma)\sigma_{G}(\xi)^{-1}\in \hat{G}'_{1}g_{1}(\sigma)$. Il en r\'esulte que les images r\'eciproques de $W_{k_{0}}$ dans ${\cal G}'_{1}$ et dans ${\cal G}'_{1,\xi}$ sont \'egales. L'assertion (5) r\'esulte alors de la d\'efinition de $\underline{\zeta}_{\xi}$. 
    
     D'apr\`es (4) et (5), l'image de $\underline{\zeta}$ dans  $H^1(\Gamma_{k_{0}},\mathfrak{Z})$ est localement triviale en presque toute place de $k_{0}$. Puisque $\Gamma_{k_{0}}$ agit trivialement sur $\mathfrak{Z}$, les \'el\'ements de $H^1(\Gamma_{k_{0}},\mathfrak{Z})$ sont simplement des homomorphismes de $\Gamma_{k_{0}}$ dans $\mathfrak{Z}$ et un homomorphisme localement trivial en presque toute place est trivial. Donc l'image de $\underline{\zeta}$ dans $H^1(\Gamma_{k_{0}},\mathfrak{Z})$ est triviale.  On a la suite exacte 
     $$1\to H^1(\Gamma_{k_{0}/k},\mathfrak{Z})\to H^1(\Gamma_{k},\mathfrak{Z})\to H^1(\Gamma_{k_{0}},\mathfrak{Z})$$
     donc $\underline{\zeta}\in H^1(\Gamma_{k_{0}/k},\mathfrak{Z})$. De m\^eme $\underline{\zeta}_{\xi}\in H^1(\Gamma_{k_{0}/k},\mathfrak{Z})$ pour tout $\xi\in \Xi$. Puisque $\Gamma_{k_{0}/k }$ est cyclique, il existe $v\in V'_{k}$ tel que $v$ soit inerte dans $k_{0}$ et $\Gamma_{k_{0,v}/k_{v}}=\Gamma_{k_{0}/k}$. L'\'egalit\'e $\underline{\zeta}_{v}=\underline{\zeta}_{\xi_{v},v}$ entra\^{\i}ne alors $\underline{\zeta}=\underline{\zeta}_{\xi_{v}}$. Cela d\'emontre (3) et le th\'eor\`eme. $\square$    
     
     \section{Exemples}\label{exemples}
   {\bf Exemple 1.} 
   
    Consid\'erons le corps   $k={\mathbb Q}_{p}$ pour un nombre premier $p\geq3$ tel que $-1$ ne soit pas un carr\'e dans $k$. Soit $\hat{G}=SL(4,{\mathbb C})$,  muni de l'action triviale de $\Gamma_{k}$. On d\'efinit  l'automorphisme habituel $\theta$ de $\hat{G}$ par $\theta(g)=J{^tg}^{-1}J^{-1}$, o\`u $J$ est la matrice antidiagonale de coefficients altern\'es $1$ et $-1$ (les coefficients $J_{m,n}$ sont nuls pour $m+n\not=5$ et valent $(-1)^{m}$ pour $m+n=5$). Soit $s$ l'\'el\'ement diagonal dont les coefficients $s_{m,m}$ valent $i$ (la racine carr\'ee de $-1$ dans ${\mathbb C}$) pour $m=1,2$ et $-i$ pour $m=3,4$.  Notons $u_{ad}$ l'image dans $\hat{G}_{AD}=PGL(4,{\mathbb C})$ de la sym\'etrie \'el\'ementaire ${\bf u}\in GL(4,{\mathbb C})$ qui \'echange les deux vecteurs de base standard centraux: ${\bf u}_{m,n}=1$  si $m=n=1,4$ ou $(m,n)=(2,3)$ ou $(m,n)=(3,2)$, les autres coefficients sont nuls. Soit $k'$ une extension quadratique ramifi\'ee de $k$.  Posons $\hat{G}'=Z_{\hat{G}}(s\theta)^0$ et notons ${\cal G}'_{ad}$ l'ensemble des $(g_{ad},w)\in {^LG}_{AD}$ tels que $g_{ad}\in \hat{G}'_{ad}$ si $w\in W_{k'}$ et $g_{ad}\in \hat{G}'_{ad}u_{ad}$ si $w\in W_{k}-W_{k'}$. Le couple $({\cal G}'_{ad},s_{ad}\theta)$ est une donn\'ee pr\'e-endoscopique pour $(\hat{G}_{AD},\theta)$. Montrons que
     
     (1) cette pr\'e-donn\'ee ne se rel\`eve pas en une pr\'e-donn\'ee pour $(\hat{G},\theta)$. 
     
     Le groupe $\hat{T}^{\theta}$ est connexe. On a $Z(\hat{G})\simeq \mu_{4}({\mathbb C})$ et $\hat{T}^{\theta}\cap Z(\hat{G})\simeq \mu_{2}({\mathbb C})$.Le groupe $\mathfrak{Z}$  du paragraphe  \ref{presquesimple}  est $\mu_{4}({\mathbb C})/\mu_{2}({\mathbb C})\simeq \{\pm 1\}$. Fixons une racine primitive $8$-i\`eme de l'unit\'e dans ${\mathbb C}$ que l'on note $\rho$. Notons $u$ le produit de ${\bf u}$ et de la matrice diagonale de coefficients $\rho$. Alors $u$ est un rel\`evement de $u_{ad}$ dans $\hat{G}$. Supposons que $({\cal G}'_{ad},s_{ad}\theta)$ se rel\`eve en une pr\'e-donn\'ee $({\cal G}',s\theta)$ de $(\hat{G},\theta)$. Alors il existe une application continue $z:W_{k}\to Z(\hat{G})$ telle que 
     $${\cal G}'=\{(z(w)h,w); w\in W_{k'}, h\in \hat{G}'\}\cup \{(z(w)hu; w\in W_{k}-W_{k'}, h\in \hat{G}'\}.$$
     Notons $\zeta$ l'application compos\'ee $W_{k}\stackrel{z}{\to }Z(\hat{G})\to \mathfrak{Z}=\{\pm 1\}$. L'application $\zeta$ est bien d\'etermin\'ee. Puisque ${\cal G}'$ est un groupe et que $u^2$ est la matrice diagonale de coefficients $\rho^2$, laquelle s'envoie sur $-1\in \mathfrak{Z}$ on doit avoir
     
    (2)  $\zeta(ww')=\zeta(w'w)=\zeta(w')\zeta(w)$ pour tous $w\in W_{k}$ et $w'\in W_{k'}$;
     
     (3) $\zeta(w^2)=-\zeta(w)^2=-1$ pour tout $w\in W_{k}-W_{k'}$. 
     
      En particulier, la restriction de $\zeta$ \`a $W_{k'}$ est un caract\`ere et la fonction $\zeta$ tout enti\`ere se factorise par le noyau de celui-ci. Donc $\zeta$ se factorise par le groupe de Weyl relatif $W_{k'/k}$. Notons $N_{k'/k}:k^{'\times}\to k^{\times}$ la norme. Fixons un \'el\'ement $a\in k^{\times}-N_{k'/k}(k^{'\times})$.  On sait qu'il y a une suite exacte
      $$1\to k^{'\times} \to W_{k'/k}\to \Gamma_{k'/k}\to 1$$
      et que l'on peut relever l'\'el\'ement non trivial de $\Gamma_{k'/k}$ en un \'el\'ement $w_{0}\in W_{k'/k}$ tel que $w_{0}^2=a$. La relation (2) implique non seulement que la restriction de $\zeta$ \`a $k^{'\times}$ est un caract\`ere, mais aussi que celui-ci est invariant par l'action galoisienne de $\Gamma_{k'/k}$. Il se factorise donc par la norme, c'est-\`a-dire qu'il existe un caract\`ere quadratique $\chi$ du groupe $N_{k'/k}(k^{'\times})$ tel que $\zeta(w)=\chi\circ N_{k'/k}(w)$ pour $w\in k^{'\times}$. La relation (3) appliqu\'ee \`a $w_{0}$ dit que $\zeta(a)=-1$, donc $\chi(a^2)=\chi\circ N_{k'/k}(a)=-1$. Mais nos hypoth\`eses impliquent que  $-1\not\in N_{k'/k}(k^{'\times})$. On peut donc choisir $a=-1$ et on obtient $\chi(1)=-1$, ce qui est impossible. Cela prouve (1).

     \bigskip
     
     {\bf Exemple 2.} 
     
     Consid\'erons maintenant le corps $k={\mathbb Q}$. 
  D'apr\`es \cite{S} exemple 5.6 page 35, on peut trouver un tore $T$ d\'efini sur $k$ tel que le noyau $ker^1(\Gamma_{k},T)$ de l'homomorphisme
     $$H^1(\Gamma_{k},T)\to \prod_{v\in V_{k}}H^1(\Gamma_{k_{v}},T)$$
     ait deux \'el\'ements. Notons ${\mathbb  A}$ l'anneau des ad\`eles de $k$.  On a la suite exacte
$$T({\mathbb A})\to H^0( {\mathbb A}/k,T)\to ker^1(\Gamma_{k},T)\to \{1\},$$
avec les notations de \cite{KS} appendice D.1. 
Donc $H^0({\mathbb A}/k,T)/im(T({\mathbb A}))$ a aussi deux \'el\'ements. Il en est de m\^eme du noyau de l'homomorphisme 
$$Hom_{cont}(H^0({\mathbb  A}/k,T),{\mathbb C}^{\times})\to Hom_{cont}(T({\mathbb  A}),{\mathbb C}^{\times}).$$
D'apr\`es \cite{KS} lemme D.2.A,  cela entra\^{\i}ne que le noyau l'homomorphisme
$$H^1(W_{k},\hat{T})\to \prod_{v\in V_{k}}H^1(W_{k_{v}},\hat{T})$$
  a aussi  deux \'el\'ements. 
        Notons $t_{1}$ et $t_{2}$ des cocycles repr\'esentant les deux \'el\'ements du noyau. Munissons $\hat{T}$ de l'automorphisme $\theta$ d\'efini par $\theta(t)=t^{-1}$. Posons $\tilde{s}=\theta$ et d\'efinissons les deux sous-groupes ${\cal G}'_{i}=\{(t_{i}(w),w);w\in W_{k}\}$ de $^LT$ pour $i=1,2$. On montre ais\'ement 
  que les couples $({\cal G}'_{1},\tilde{s})$ et $({\cal G}'_{2},\tilde{s})$ sont deux pr\'e-donn\'ees endoscopiques elliptiques pour $(\hat{T},\theta)$ qui sont localement \'equivalentes en toute place de $k$ mais ne sont pas \'equivalentes.

  \bigskip
  
  CNRS Institut de Math\'ematiques de Jussieu-Paris Rive Gauche
  
  jean-loup.waldspurger@imj-prg.fr
          
 \end{document}